\documentclass[reqno]{amsart}
\usepackage{setspace}
\usepackage{graphicx}
\usepackage{amsmath}
\usepackage{amsfonts}
\usepackage{amssymb}
\begin{document}

\centerline{\huge \bf On Generalized Kissing Numbers of}
\centerline{\huge \bf Convex Bodies (\uppercase\expandafter{\romannumeral2})}

\bigskip
\centerline{\large Yiming Li and Chuanming Zong}

{\large
\vspace{0.8cm}
\centerline{\begin{minipage}{12.8cm}
{\bf Abstract.} In 1694, Gregory and Newton discussed the problem to determine the kissing number of a rigid material ball. This problem and its higher dimensional generalization have been studied by many mathematicians, including Minkowski, van der Waerden, Hadwiger, Swinnerton-Dyer, Watson, Levenshtein, Odlyzko, Sloane and Musin. Recently, Li and Zong introduced and studied the generalized kissing numbers of convex bodies. As a continuation of this project, in this paper we obtain the exact generalized kissing numbers $\kappa_\alpha^*(B^n)$ of the $n$-dimensional balls for $3\le n\le 8$ and $\alpha =2\sqrt{3}-2$. Furthermore, the lattice kissing number of a four-dimensional cross-polytope is determined.
\end{minipage}}



\vspace{1cm}
\noindent
{\LARGE\bf 1. Introduction}

\bigskip
\noindent
In 1694, Gregory and Newton discussed the following problem: {\it Can a rigid material ball be brought into contact with thirteen other such balls of the same size}? Gregory believed \lq\lq yes", while Newton thought \lq\lq no". The solution of this problem has a complicated history! Several authors claimed proofs that {\it the largest number of nonoverlapping unit balls which can be brought into contact with a fixed one is twelve} (see Hoppe \cite{Ho01}, G\"{u}nter \cite{Gu01}, Sch\"{u}tte and van der Waerden \cite{Sc01} and Leech \cite{Le02}). However, only Sch\"{u}tte and van der Waerden's proof is complete!

Let $K$ be an $n$-dimensional convex body and let $B^n$ denote the $n$-dimensional unit ball centered at the origin ${\bf o}$ of $\mathbb{E}^n$. Let $\kappa(K)$ and $\kappa^*(K)$ denote the translative kissing number and the lattice kissing number of $K$, respectively. In other words, $\kappa(K)$ is the maximum number of nonoverlapping translates $K+{\bf x}$ that can touch $K$ at its boundary, and $\kappa^*(K)$ is defined similarly, with the restriction that the translates are members of a lattice packing of $K$. The kissing numbers $\kappa(K)$ and $\kappa^*(K)$, especially for $\kappa(B^n)$ and $\kappa^*(B^n)$, have been studied by many authors (see \cite{Ba01, Le01, Mu01, Mu02, Od01, Wa01}).

Define $\kappa^\star(K)$ to be the second lattice kissing number of $K$. In other words, $\kappa^\star(K)$ is the second largest number $m$ such that there is a lattice packing $K+\Lambda $ in which $K$ touches $m$ others. It is easy to see that
$$\kappa^\star(K)\leq\kappa^*(K)-2.$$

There is a remarkable relation between lattice packings of $B^n$ and positive definite quadratic forms $F_n$ of $n$ variables. In fact, by studying the numbers of minimum integer points of $F_n$, Watson \cite{Wa01} obtained the following results:

\medskip
\centerline{
\renewcommand\arraystretch{1.5}
\begin{tabular}{|c|c|c|c|c|c|c|c|c|c|c|}
\hline $n$ & $2$ & $3$ & $4$ & $5$ & $6$ & $7$ & $8$ & $9$\\
\hline $\kappa^*(B^n)$ & $6$ & $12$ & $24$ & $40$ & $72$ & $126$ & $240$ & $272$\\
\hline $\kappa^\star(B^n)$ & $4$ & $10$ & $20$ & $32$ & $60$ & $92$ & $150$ & $??$\\
\hline
\end{tabular}}

\medskip
Recently, Li and Zong \cite{Li01} defined $\kappa^*_\alpha(C)$ to be the maximum number of translates $C+{\bf x}$ which can be packed into the region $(3+\alpha)C \setminus {\rm int}(C)$ and the translates are members of a lattice packing of $C$, where $C$ is a centrally symmetric convex body centered at ${\bf o}$ and ${\rm int}(C)$ denotes the interior of $C$. Among other things, they proved that, for $\alpha=2\sqrt2-2$, we have

\medskip
\centerline{
\renewcommand\arraystretch{1.5}
\begin{tabular}{|c|c|c|c|c|c|c|c|c|c|c|}
\hline $n$ & $2$ & $3$ & $4$ & $5$ & $6$ & $7$ & $8$\\
\hline $\kappa^*_\alpha(B^n)$ & $8$ & $20$ & $50$ & $??$ & $??$ & $??$ & $2400$\\
\hline
\end{tabular}}

\medskip
\noindent
It is interesting to notice that, the optimal lattices for $\kappa^*_\alpha(B^8)$ and $\kappa^*(B^8)$ are the same. However, this is not true in 
$\mathbb{E}^2$, $\mathbb{E}^3$ and $\mathbb{E}^4$.

In this paper, we prove the following results:

\medskip
\noindent
{\bf Theorem 1.} {\it If $\alpha=2\sqrt3-2$, we have

\medskip
\centerline{
\renewcommand\arraystretch{1.5}
\begin{tabular}{|c|c|c|c|c|c|c|c|c|c|c|}
\hline $n$ & $2$ & $3$ & $4$ & $5$ & $6$ & $7$ & $8$\\
\hline $\kappa^*_\alpha(B^n)$ & $12$ & $42$ & $144$ & $370$ & $1062$ & $2954$ & $9120$\\
\hline
\end{tabular}}

\medskip
\noindent
where the optimal lattices are the corresponding optimal lattices of $\kappa^*(B^n)$, respectively.}

\medskip
\noindent
{\bf Theorem 2.} {\it Let $C_3$ denote a three-dimensional cross-polytope, then we have
$$\kappa^*(C_3)=18,\quad \kappa^\star(C_3)=14.$$}

\medskip
\noindent
{\bf Theorem 3.} {\it Let $C_4$ denote a four-dimensional cross-polytope, then we have
$$\kappa^*(C_4)=40.$$}

\medskip
\noindent
{\bf Remark 1.} Let $\ell (\Lambda)$ denote the length of the shortest non-zero vector of a given $n$-dimensional lattice $\Lambda$ and let $\gamma (n)$ be a function of $n$. To find a non-zero lattice vector ${\bf v}$ satisfying $\| {\bf v}\|\le \gamma (n) \ell (\Lambda)$ is a fundamental problem in lattice-based cryptography (see \cite{Micc02}). It is easy to see that the generalized kissing numbers provide estimations for the number of such vectors.

\medskip
\noindent
{\bf Remark 2.} Similar to the tetrahedral case (see \cite{Zo02, Zo03}),  the density of the lattice packing $C_4+\Lambda$ attaining $\kappa^*(C_4)=40$ is only $2/3$ while the density of some lattice packing $C_4+\Lambda$ in which every one touches 26 others is $32/45$.

\vspace{1cm}
\noindent
{\LARGE\bf 2. Proof of Theorem 1}

\bigskip
\noindent
It is well-known that each centrally symmetric convex body $C$ centered at ${\bf o}$ defines a metric $\|\cdot \|_C$ on $\mathbb{R}^n$ by
$$\|{\bf x}, {\bf y}\|_C=\|{\bf x}-{\bf y}\|_C=\min \{ r:\ r>0,\ {\bf x}-{\bf y}\in rC\}.$$
Especially, we use $||\cdot||$ to denote the metric defined by $B^n$. For an non-negative number $\alpha$ and a packing lattice $\Lambda $ of $B^n$, we define
$$X(\alpha, \Lambda)=\{{\bf v}:\ 2\leq ||{\bf v}||\leq2+\alpha,\ {\bf v}\in\Lambda\}.$$

In \cite{Li01}, Li and Zong proved the following result.

\medskip
\noindent
{\bf Lemma 2.1.} {\it If $\alpha<2\sqrt3-2$, the set $X(\alpha, \Lambda)$ contains no four collinear points. If $\alpha=2\sqrt3-2$, the set $X(\alpha, \Lambda)$ has four collinear points ${\bf v}_1, {\bf v}_2, {\bf v}_3$ and ${\bf v}_4$ if and only if
$${\bf v}_1=(\sqrt3, 3),\quad {\bf v}_2=(\sqrt3, 1),\quad {\bf v}_3=(\sqrt3, -1)\quad {\rm and} \quad {\bf v}_4=(\sqrt3, -3),$$
up to rotation, reflection and re-enumeration.}

\medskip
For two vectors ${\bf v}$ and ${\bf v}'$ which belongs to an $n$-dimensional lattice $\Lambda$, we say ${\bf v}$ are modulo 3 equivalent to ${\bf v}'$, if and only if
$$\mbox{$\frac{1}{3}$}({\bf v}-{\bf v}')\in\Lambda.$$
In other words, both four collinear points ${\bf v}, \frac{2}{3}{\bf v}+\frac{1}{3}{\bf v}', \frac{2}{3}{\bf v}'+\frac{1}{3}{\bf v}, {\bf v}'$ belongs to $\Lambda$. Clearly, this relation divides the points of $\Lambda$ into $3^n$ classes. 

Next we prove the following lemmas.

\medskip
\noindent
{\bf Lemma 2.2.} {\it If $\alpha<2\sqrt3-2$, we have
$$\kappa^*_\alpha(B^n)\leq3^n-1.$$}

\begin{proof}
On the contrary, suppose that there are such a positive $\alpha$ which less than $2\sqrt3-2$ and a suitable lattice $\Lambda$ satisfying
$${\rm card}\{X(\alpha, \Lambda)\}\geq 3^n+1.$$

Since ${\rm card}\{X(\alpha, \Lambda)\}\geq 3^n+1,$ two of them must belong to the same modulo 3 equivalent class of $\Lambda$. Without loss of generality, we may assume ${\bf v}, {\bf v}'\in X(\alpha, \Lambda)$ are modulo 3 equivalent. Then, we have four collinear points
$${\bf v},\quad \mbox{$\frac{2}{3}$}{\bf v}+\mbox{$\frac{1}{3}$}{\bf v}',\quad \mbox{$\frac{2}{3}$}{\bf v}'+\mbox{$\frac{1}{3}$}{\bf v},\quad {\bf v}'\ \in X(\alpha, \Lambda),$$
which contradicts Lemma 2.1. Therefore, for $\alpha<2\sqrt3-2$, we have
$$\kappa^*_\alpha(B^n)\leq3^n-1.$$
Lemma 2.2 is proved.
\end{proof}

\medskip
\noindent
{\bf Lemma 2.3.} {\it If $\alpha=2\sqrt3-2$, a modulo 3 equivalent class of $\Lambda$ can contain exactly 0, 1 or 3 points of $X(\alpha, \Lambda)$.}

\begin{proof}
Suppose there exist ${\bf v}_1, {\bf v}_4\in X(\alpha, \Lambda)$ which belong to the same modulo 3 equivalent class of $\Lambda$. By Lemma 2.1, we may assume
$${\bf v}_1=(\sqrt3, 3, 0,\ldots,0),\quad {\bf v}_4=(\sqrt3, -3, 0,\ldots,0)\ \in\Lambda$$
and
$${\bf v}_2=(\sqrt3, 1, 0,\ldots,0),\quad {\bf v}_3=(\sqrt3, -1, 0,\ldots,0)\ \in\Lambda,$$
without loss of generality.

Define
$$\pi_0:\ \{{\bf v}=(v_1, v_2,\ldots,v_n):\ v_3=v_4=\ldots=v_n=0\}.$$
It is easy to see that the two-dimensional lattice $\Lambda'=\pi_0\cap\Lambda$ is generated by $\{{\bf v}_2, {\bf v}_3\}$ and $-({\bf v}_2+{\bf v}_3)$ are modulo 3 equivalent to ${\bf v}_1$ and ${\bf v}_4$.

Suppose there exist another point ${\bf v}=(v_1, v_2,\ldots,v_n)\in X(\alpha, \Lambda)$ which are modulo 3 equivalent to ${\bf v}_1, {\bf v}_4$ and $-({\bf v}_2+{\bf v}_3)$. By Lemma 2.1, we have
$$||{\bf v}||=2\sqrt3,\quad ||{\bf v}, {\bf v}_1||=||{\bf v}, {\bf v}_4||=||{\bf v}, -({\bf v}_2+{\bf v}_3)||=6,$$
namely
$$\left\{
\begin{aligned}
& v^2_1+v^2_2+\ldots+v^2_n=12,\\
& (v_1-\sqrt3)^2+(v_2-3)^2+v_3^2+\ldots+v^2_n=36,\\
& (v_1-\sqrt3)^2+(v_2+3)^2+v_3^2+\ldots+v^2_n=36,\\
& (v_1+2\sqrt3)^2+v_2^2+v_3^2+\ldots+v^2_n=36,\\
\end{aligned}
\right.
$$
which do not exist solution. Therefore, when $\alpha=2\sqrt3-2$, a modulo 3 equivalent class of $\Lambda$ can contain exactly 0, 1 or 3 points of $X(\alpha, \Lambda)$.

Lemma 2.3 is proved.
\end{proof}

\medskip
For a packing lattice $\Lambda$ of $B^n$, define $\kappa(B^n,\Lambda)={\rm card}\{{\bf v}:\ ||{\bf v}||=2,\ {\bf v}\in\Lambda\}$. Then, we have:

\medskip
\noindent
{\bf Lemma 2.4.} {\it If $\alpha=2\sqrt3-2$, we have
$${\rm card}\{X(\alpha, \Lambda)\}\leq3^n-1+\mbox{$\frac{1}{9}$}(\kappa(B^n, \Lambda))^2.$$}

\begin{proof}
For $\alpha=2\sqrt3-2$ and a packing lattice $\Lambda$ of $B^n$, denote the numbers of modulo 3 equivalent classes of $\Lambda$ which contain exactly $i$ points of $X(\alpha, \Lambda)$ by $m_i$. Obviously, the modulo 3 equivalent class which contain ${\bf o}$ cannot contain the points of $X(\alpha, \Lambda)$. Therefore, we have
$$m_3+m_1\leq3^n-1.\eqno(2.1)$$
By Lemma 2.3, we have
$${\rm card}\{X(\alpha, \Lambda)\}=3m_3+m_1.\eqno(2.2)$$

Define two collections of sets
\begin{eqnarray*}
C(X(\alpha, \Lambda))=\Big\{\hspace{-0.4cm} &\big\{ \hspace{-0.35cm}&{\bf v}_1,\ {\bf v}_2,\ {\bf v}_3,\ \mbox{$\frac{2}{3}$}{\bf v}_{1}+\mbox{$\frac{1}{3}$}{\bf v}_{2},\ \mbox{$\frac{2}{3}$}{\bf v}_{2}+\mbox{$\frac{1}{3}$}{\bf v}_{1},\ \mbox{$\frac{2}{3}$}{\bf v}_{1}+\mbox{$\frac{1}{3}$}{\bf v}_{3},\ \mbox{$\frac{2}{3}$}{\bf v}_{3}+\mbox{$\frac{1}{3}$}{\bf v}_{1},\\ 
\hspace{-0.1cm} &\mbox{$\frac{2}{3}$} \hspace{-0.35cm}& {\bf v}_{2}+\mbox{$\frac{1}{3}$}{\bf v}_{3},\ \mbox{$\frac{2}{3}$}{\bf v}_{3}+\mbox{$\frac{1}{3}$}{\bf v}_{2}\}:\ {\bf v}_1, {\bf v}_2, {\bf v}_3\in X(\alpha,\Lambda)\ {\rm are\ modulo\ 3\ equivalent}\Big\}
\end{eqnarray*}
and
$$SV_2(\Lambda)=\big\{\{{\bf v}, {\bf v}'\}:\ ||{\bf v}||=||{\bf v}'||=2,\ ||{\bf v},{\bf v}'||=2,\ {\bf v}, {\bf v}'\in\Lambda\big\}.$$
Clearly, we have
$${\rm card}\{C(X(\alpha, \Lambda))\}=m_3. \eqno(2.3)$$

By Lemma 2.1 and Lemma 2.3, a set $V\in C(X(\alpha, \Lambda))$, if and only if
$$||\mbox{$\frac{2}{3}$}{\bf v}_{i_1}+\mbox{$\frac{1}{3}$}{\bf v}_{i_2}||=2$$
holds for $\{i_1,i_2\}\subset\{1,2,3\}$. In other words,
$$||\mbox{$\frac{2}{3}$}{\bf v}_{i_1}+\mbox{$\frac{1}{3}$}{\bf v}_{i_2},\ {\bf o}||=||\mbox{$\frac{2}{3}$}{\bf v}_{i_1}+\mbox{$\frac{1}{3}$}{\bf v}_{i_2},\ \mbox{$\frac{2}{3}$}{\bf v}_{i_1}+\mbox{$\frac{1}{3}$}{\bf v}_{i_3}||=||\mbox{$\frac{2}{3}$}{\bf v}_{i_1}+\mbox{$\frac{1}{3}$}{\bf v}_{i_2},\ \mbox{$\frac{2}{3}$}{\bf v}_{i_2}+\mbox{$\frac{1}{3}$}{\bf v}_{i_1}||=2$$
holds for $\{i_1, i_2, i_3\}=\{1, 2, 3\}$, see Figure 1.

\begin{figure}[h]
\centering
\includegraphics[height=6.5cm]{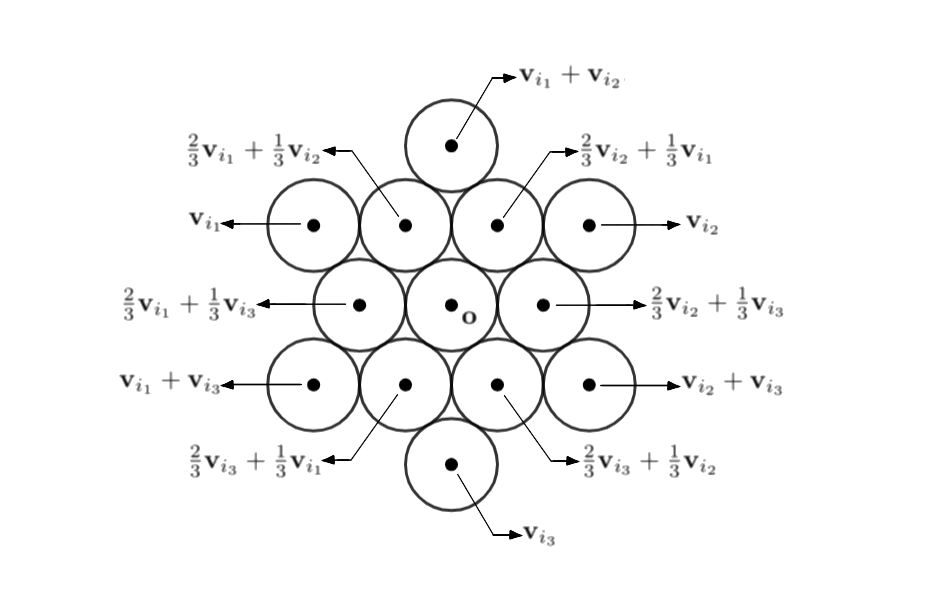}
\caption{$||\mbox{$\frac{2}{3}$}{\bf v}_{i_1}+\mbox{$\frac{1}{3}$}{\bf v}_{i_2}||=2$ holds for $\{i_1,i_2\}\subset\{1,2,3\}$.}
\label{fig1}
\end{figure}

Since a set which belong to $C(X(\alpha, \Lambda))$ contain exactly six sets of $SV_2(\Lambda)$, a set of $SV_2(\Lambda)$ contained in exactly two sets of $C(X(\alpha, \Lambda))$ (see Figure 1), we have
$${\rm card}\{C(X(\alpha, \Lambda))\}=\mbox{$\frac{1}{3}$}{\rm card}\{SV_2(\Lambda)\}. \eqno(2.4)$$

For a lattice point ${\bf v}\in\Lambda$ which length is 2, define
$$S({\bf v},\Lambda)=\{{\bf v}':\ ||{\bf v}', {\bf v}||=2,\ ||{\bf v}'||=2,\ {\bf v}'\in\Lambda\}$$
and
$$t({\bf v},\Lambda)={\rm card}\{S({\bf v},\Lambda)\}.$$
Clearly, we have
$${\rm card}\{SV_2(\Lambda)\}=\frac{1}{2}\sum_{{\bf v}\in\Lambda,\ ||{\bf v}||=2}t({\bf v}, \Lambda). \eqno(2.5)$$
Combining with (2.1)-(2.4), we have
\begin{eqnarray*}
{\rm card}\{X(\alpha, \Lambda)\}\hspace{-0.2cm} &\leq \hspace{-0.2cm}&3^n-1+2m_3\\
\hspace{-0.2cm} &= \hspace{-0.2cm}&3^n-1+\frac{1}{3}\sum_{{\bf v}\in\Lambda,\ ||{\bf v}||=2}t({\bf v}, \Lambda).
\end{eqnarray*}
$$\eqno(2.6)$$

For a packing lattice $\Lambda$ of $B^n$, define
$$t_\Lambda=\max_{{\bf v}\in\Lambda,\ ||{\bf v}||=2}t({\bf v},\Lambda).\eqno(2.7)$$
For convenience, we suppose lattice point ${\bf v}_0\in\Lambda$ which length is 2, satisfy $t({\bf v}_0,\Lambda)=t_\Lambda$. Denote $S({\bf v}_0,\Lambda)=\{{\bf v}_1, {\bf v}_0-{\bf v}_1,\ldots, {\bf v}_{t_\Lambda/2}, {\bf v}_0-{\bf v}_{t_\Lambda/2}\}$. Without loss of generality, assume
$$t({\bf v}_1,\Lambda)=\max_{{\bf v}\in S({\bf v}_0,\Lambda)}t({\bf v}, \Lambda).$$
For convenience, we denote $t({\bf v}_1,\Lambda)$ by $t'_\Lambda$. Denote
$$S({\bf v}_1,\Lambda)=\{{\bf v}_0,\ {\bf v}_1-{\bf v}_0,\ {\bf v}_1',\ {\bf v}_1-{\bf v}_1',\ {\bf v}_2',\ {\bf v}_1-{\bf v}_2',\ldots, {\bf v}'_{(t_\Lambda'-2)/2},\ {\bf v}_1-{\bf v}_{(t_\Lambda'-2)/2}'\}.$$
It is easy to see that, ${\bf v}_i'$ and ${\bf v}_1-{\bf v}_i'$ cannot belong to $S({\bf v}_0,\Lambda)\cup \big(-S({\bf v}_0,\Lambda)\big)\cup\{\pm{\bf v}_0\}$ simultaneously, $i\in\{1,\ldots,(t_{\Lambda}'-2)/2\}$. Therefore, we have
$$\kappa(B^n,\Lambda)\geq2+2t_\Lambda+2(\mbox{$\frac{t'_\Lambda}{2}$}-1)=2t_\Lambda+t_{\Lambda}',$$
which means
\begin{eqnarray*}
\sum_{{\bf v}\in\Lambda,\ ||{\bf v}||=2}t({\bf v}, \Lambda)\hspace{-0.2cm} &= \hspace{-0.2cm}&\sum_{{\bf v}\notin S({\bf v}_0,\Lambda)\cup(-S({\bf v}_0,\Lambda))}t({\bf v}, \Lambda)+\sum_{{\bf v}\in S({\bf v}_0,\Lambda)\cup(-S({\bf v}_0,\Lambda))}t({\bf v}, \Lambda)\\
\hspace{-0.2cm} &\leq \hspace{-0.2cm}&(\kappa(B^n,\Lambda)-2t_\Lambda)t_\Lambda+2t_\Lambda t'_\Lambda\\
\hspace{-0.2cm} &\leq \hspace{-0.2cm}&(\kappa(B^n,\Lambda)-2t_\Lambda)t_\Lambda+2t_\Lambda(\kappa(B^n,\Lambda)-2t_\Lambda)\\
\hspace{-0.2cm} &= \hspace{-0.2cm}&-6t_\Lambda^2+3\kappa(B^n,\Lambda)t_\Lambda.
\end{eqnarray*}
$$\eqno(2.8)$$

If $t_\Lambda<\frac{1}{3}\kappa(B^n,\Lambda)$, by (2.6) and (2.7), we have
$${\rm card}\{X(\alpha, \Lambda)\}\leq3^n-1+\mbox{$\frac{1}{3}$}\kappa(B^n,\Lambda)t_\Lambda<3^n-1+\mbox{$\frac{1}{9}$}(\kappa(B^n,\Lambda))^2. \eqno(2.9)$$
On the other hand, if $t_\Lambda\geq\frac{1}{3}\kappa(B^n,\Lambda)$, by (2.6) and (2.8), we have
\begin{eqnarray*}
{\rm card}\{X(\alpha, \Lambda)\}\hspace{-0.2cm} &\leq \hspace{-0.2cm}&3^n-1+(-2t_\Lambda^2+\kappa(B^n,\Lambda)t_\Lambda)\\
\hspace{-0.2cm} &\leq \hspace{-0.2cm}&3^n-1+(-2)\cdot(\mbox{$\frac{1}{3}$}\kappa(B^n,\Lambda))^2+\kappa(B^n,\Lambda)\cdot\mbox{$\frac{1}{3}$}\kappa(B^n,\Lambda)\\
\hspace{-0.2cm} &= \hspace{-0.2cm}&3^n-1+\mbox{$\frac{1}{9}$}(\kappa(B^n,\Lambda))^2.
\end{eqnarray*}
Lemma 2.4 is proved.
\end{proof}

\medskip
\noindent
{\bf Remark 2.1.}
{\rm Suppose $\Lambda_0$ is the unique lattice which satisfy
$$\kappa(B^n,\Lambda_0)=\kappa^*(B^n),$$
up to rotation and reflection. If for $\alpha=2\sqrt3-2$, we have
$${\rm card}\{X(\alpha, \Lambda_0)\}=3^n-1+\mbox{$\frac{1}{9}$}(\kappa^*(B^n))^2,$$
then by Lemma 2.4, we have
$$\kappa^*_\alpha(B^n)={\rm card}\{X(\alpha, \Lambda_0)\}=3^n-1+\mbox{$\frac{1}{9}$}(\kappa^*(B^n))^2,$$
and it can be attained if and only if the corresponding lattice is $\Lambda_0$, up to rotation and reflection.

However, if
$${\rm card}\{X(\alpha, \Lambda_0)\}<3^n-1+\mbox{$\frac{1}{9}$}(\kappa^*(B^n))^2,$$
by Lemma 2.4 and the definition of $\kappa^\star(B^n)$, we have
\begin{eqnarray*}
\kappa^*_\alpha(B^n)\hspace{-0.2cm} &= \hspace{-0.2cm}&\max\Big\{{\rm card}\{X(\alpha, \Lambda_0)\},\ \max_{\Lambda\neq\Lambda_0}{\rm card}\{X(\alpha, \Lambda)\}\Big\}\\
\hspace{-0.2cm} &\leq \hspace{-0.2cm}&\max\Big\{{\rm card}\{X(\alpha, \Lambda_0)\},\ 3^n-1+\mbox{$\frac{1}{9}$}(\kappa^\star(B^n))^2\Big\}.
\end{eqnarray*}
In other words, if ${\rm card}\{X(\alpha, \Lambda_0)\}> 3^n-1+\mbox{$\frac{1}{9}$}(\kappa^\star(B^n))^2$, we still can determine that
$$\kappa^*_\alpha(B^n)={\rm card}\{X(\alpha, \Lambda_0)\},$$
and it can be attained if and only if the corresponding lattice is $\Lambda_0$, up to rotation and reflection.

Thus, for some problems, the concept of $\kappa^\star(B^n)$ will be useful, especially when the gap between $\kappa^{*}(B^n)$ and $\kappa^\star(B^n)$ is large.}

\medskip
\noindent
{\bf Theorem 1.} {\it If $\alpha=2\sqrt3-2$, we have

\medskip
\centerline{
\renewcommand\arraystretch{1.5}
\begin{tabular}{|c|c|c|c|c|c|c|c|c|c|c|}
\hline $n$ & $2$ & $3$ & $4$ & $5$ & $6$ & $7$ & $8$\\
\hline $\kappa^*_\alpha(B^n)$ & $12$ & $42$ & $144$ & $370$ & $1062$ & $2954$ & $9120$\\
\hline
\end{tabular}}

\medskip
\noindent
where the optimal lattices are the corresponding optimal lattices of $\kappa^*(B^n)$, respectively.}
\begin{proof}
As usual, we write
$$A_n=\Big\{(v_0,v_1,v_2,\ldots,v_n):\ v_i\in Z;\ \sum v_i=0\Big\}$$
and
$$D_n=\Big\{(v_1,v_2,\ldots,v_n):\ v_i\in Z;\ \sum v_i\in 2Z\Big\}.$$

Fix a lattice $D_5,\ D_6,\ D_7$ in the hyperplane $H_6=\{{\bf v}\in\mathbb{E}^6:\ v_6=0\}$, $H_7=\{{\bf v}\in\mathbb{E}^7:\ v_7=0\}$, $H_8=\{{\bf v}\in\mathbb{E}^8:\ v_8=0\}$ in $\mathbb{E}^6$, $\mathbb{E}^7$, $\mathbb{E}^8$, respectively. We write
$$E_6=\{z{\bf v}_6+D_5:\ z\in\mathbb{Z}\},\ E_7=\{z{\bf v}_7+D_6:\ z\in\mathbb{Z}\},\ E_8=\{z{\bf v}_8+D_7:\ z\in\mathbb{Z}\},$$
where ${\bf v}_6=(\frac{1}{2}, \frac{1}{2}, \frac{1}{2}, \frac{1}{2}, \frac{1}{2}, \frac{\sqrt3}{2})$, ${\bf v}_7=(\frac{1}{2}, \frac{1}{2}, \frac{1}{2}, \frac{1}{2}, \frac{1}{2}, \frac{1}{2}, \frac{\sqrt2}{2})$, ${\bf v}_8=(\frac{1}{2}, \frac{1}{2}, \frac{1}{2}, \frac{1}{2}, \frac{1}{2}, \frac{1}{2}, \frac{1}{2}, \frac{1}{2})$.

It is easy to verify that, for $\alpha=2\sqrt3-2$, we have
\begin{eqnarray*}
{\rm card}\{X(\alpha, \sqrt2A_2)\}\hspace{-0.2cm} &= \hspace{-0.2cm}&12=3^2-1+\mbox{$\frac{1}{9}$}\times(\kappa^*(B^2))^2,\\
{\rm card}\{X(\alpha, \sqrt2A_3)\}\hspace{-0.2cm} &= \hspace{-0.2cm}&42=3^3-1+\mbox{$\frac{1}{9}$}\times(\kappa^*(B^3))^2,\\
{\rm card}\{X(\alpha, \sqrt2D_4)\}\hspace{-0.2cm} &= \hspace{-0.2cm}&144=3^4-1+\mbox{$\frac{1}{9}$}\times(\kappa^*(B^4))^2.
\end{eqnarray*}
Since $\sqrt2A_2, \sqrt2A_3, \sqrt2D_4$ is the unique lattice to attain $\kappa^*(B^2), \kappa^*(B^3), \kappa^*(B^4)$ (see \cite{Zo01}), respectively, up to rotation and reflection. By Remark 2.1, we have
$$\kappa^*_\alpha(B^2)=12,\quad \kappa^*_\alpha(B^3)=42,\quad \kappa^*_\alpha(B^4)=144,$$
and it can be attained if and only if the corresponding lattice is $\sqrt2A_2, \sqrt2A_3, \sqrt2D_4$, respectively, up to rotation and reflection.

By enumeration (or see \cite{Co01}), one can verify that,
$${\rm card}\{X(\alpha, \sqrt2D_5)\}=370,\quad {\rm card}\{X(\alpha, \sqrt2E_8)\}=9120.$$
Since
\begin{eqnarray*}
3^5-1+\mbox{$\frac{1}{9}$}\times(\kappa^*(B^5))^2\hspace{-0.2cm} &= \hspace{-0.2cm}&419.77\cdots>{\rm card}\{X(\alpha, \sqrt2D_5)\},\\
3^5-1+\mbox{$\frac{1}{9}$}\times(\kappa^\star(B^5))^2\hspace{-0.2cm} &= \hspace{-0.2cm}&355.77\cdots<{\rm card}\{X(\alpha, \sqrt2D_5)\},\\
3^8-1+\mbox{$\frac{1}{9}$}\times(\kappa^*(B^8))^2\hspace{-0.2cm} &= \hspace{-0.2cm}&12960>{\rm card}\{X(\alpha, \sqrt2E_8)\},\\
3^8-1+\mbox{$\frac{1}{9}$}\times(\kappa^\star(B^8))^2\hspace{-0.2cm} &= \hspace{-0.2cm}&9060<{\rm card}\{X(\alpha, \sqrt2E_8)\},
\end{eqnarray*}
and $\sqrt2D_5, \sqrt2E_8$ is the unique lattice to attain $\kappa^*(B^5), \kappa^*(B^8)$ (see \cite{Zo01}), respectively, up to rotation and reflection. By Remark 2.1, we have
$$\kappa^*_\alpha(B^5)=370,\quad \kappa^*_\alpha(B^8)=9120,$$
and it can be attained if and only if the corresponding lattice is $\sqrt2D_5, \sqrt2E_8$, respectively, up to rotation and reflection.

In $\mathbb{E}^6$, fix a lattice $D_5$ in the hyperplane $H_6$ and take
$$D_{5,{\bf v}}=\{z{\bf v}+D_5:\ z\in\mathbb{Z}\}.$$
By the discussions of Watson (see \cite{Wa01}), for a 6-dimensional lattice $\Lambda$, $\kappa(B^6,\Lambda)\geq56$ if and only if $\Lambda=\sqrt2E_6, \sqrt2D_6$ or $\sqrt2D_{5,{\bf v}_6'}$, up to rotation and reflection, where ${\bf v}_6'=(\frac{1}{2}, \frac{1}{2}, \frac{1}{2}, \frac{1}{2}, v_5, \sqrt{1-v_5^2})$, $0\leq v_5<\frac{1}{2}$.

\medskip
By enumeration and (2.9), one can verify that:
$${\rm card}\{X(\alpha, \sqrt2E_6)\}=1062,\quad {\rm card}\{X(\alpha, \sqrt2D_6)\}=856,$$
\begin{eqnarray*}
\hspace{0.45cm}{\rm card}\{X(\alpha, \sqrt2D_{5,{\bf v}_6'})\}\hspace{-0.2cm} &\leq \hspace{-0.2cm}& 3^6-1+\mbox{$\frac{1}{3}$}\kappa(B^6,\sqrt2D_{5,{\bf v}_6'})t_{\sqrt2D_{5,{\bf v}_6'}}\\
\hspace{-0.2cm} &= \hspace{-0.2cm}&3^6-1+\mbox{$\frac{1}{3}$}\times56\times16\\
\hspace{-0.2cm} &= \hspace{-0.2cm}&1026.66\cdots,
\end{eqnarray*}
where $t_\Lambda$ for a packing lattice $\Lambda$ has already been defined, see (2.7).

Since
$${\rm card}\{X(\alpha, \sqrt2E_6)\}>\max\big\{{\rm card}\{X(\alpha, \sqrt2D_6)\},\ {\rm card}\{X(\alpha, \sqrt2D_{5,{\bf v}_6'})\},\ 3^6-1+\mbox{$\frac{1}{9}$}\times(54)^2\big\}$$
and $\sqrt2E_6$ is the unique lattice to attain $\kappa^*(B^6)$ (see \cite{Zo01}), up to rotation and reflection. By Remark 2.1, we have
$$\kappa^*_\alpha(B^6)=1062,$$
and it can be attained if and only if the corresponding lattice is $\sqrt2E_6$, up to rotation and reflection.

In $\mathbb{E}^7$, fix a lattice $E_6$ in the hyperplane $H_7$ and take
$$E_{6,{\bf v}}=\{z{\bf v}+E_6:\ z\in\mathbb{Z}\};$$
fix a lattice $D_6$ in the hyperplane $H_7$ and take
$$D_{6,{\bf v}}=\{z{\bf v}+D_6:\ z\in\mathbb{Z}\}.$$

By the discussions of Watson (see \cite{Wa01}), for a 7-dimensional lattice $\Lambda$, $\kappa(B^7,\Lambda)\geq92$ if and only if $\Lambda=\sqrt2E_7$ or $\sqrt2D_{6,{\bf v}_7'}$, up to rotation and reflection, where
$${\bf v}_7'=(\mbox{$\frac{1}{2}, \frac{1}{2}, \frac{1}{2}, \frac{1}{2}, \frac{1}{2}, v_6, \sqrt{\frac{3}{4}-{v^2_6}}$}),\quad 0\leq v_6<\mbox{$\frac{1}{2}$}.$$
Furthermore, Watson \cite{Wa01} proved that: if $\kappa(B^7,\Lambda)<92$, then we have
$$\kappa(B^7,\Lambda)\leq84,$$
and the equality holds if and only if $\Lambda=\sqrt2E_{6,{\bf v}'}$ or $\sqrt2D_7$, where ${\bf v}'$ satisfy some special conditions, up to rotation and reflection.

By enumeration and (2.9), one can verify that:
$${\rm card}\{X(\alpha, \sqrt2E_7)\}=2954,\quad {\rm card}\{X(\alpha, \sqrt2D_7)\}=1946,$$
\begin{eqnarray*}
\hspace{0.25cm}{\rm card}\{X(\alpha, \sqrt2D_{6,{\bf v}_7'})\}\hspace{-0.2cm} &\leq \hspace{-0.2cm}& 3^7-1+\mbox{$\frac{1}{3}$}\kappa(B^7,\sqrt2D_{6,{\bf v}_7'})t_{\sqrt2D_{6,{\bf v}_7'}}\\
\hspace{-0.2cm} &\leq \hspace{-0.2cm}&3^7-1+\mbox{$\frac{1}{3}$}\times92\times24\\
\hspace{-0.2cm} &= \hspace{-0.2cm}&2922,
\end{eqnarray*}
\begin{eqnarray*}
{\rm card}\{X(\alpha, \sqrt2E_{6,{\bf v}'})\}\hspace{-0.2cm} &\leq \hspace{-0.2cm}& 3^7-1+\mbox{$\frac{1}{3}$}\kappa(B^7,\sqrt2E_{6,{\bf v}'})t_{\sqrt2E_{6,{\bf v}'}}\\
\hspace{-0.2cm} &\leq \hspace{-0.2cm}&3^7-1+\mbox{$\frac{1}{3}$}\times84\times26\\
\hspace{-0.2cm} &= \hspace{-0.2cm}&2914.
\end{eqnarray*}

Since
\begin{eqnarray*}
{\rm card}\{X(\alpha, \sqrt2E_7)\}>\max\big\{\hspace{-0.35cm} &{\rm card} \hspace{-0.35cm}&\{X(\alpha, \sqrt2D_7)\},\ {\rm card}\{X(\alpha, \sqrt2D_{6,{\bf v}_7'})\},\\
\hspace{-0.35cm} &{\rm card} \hspace{-0.35cm}&\{X(\alpha, \sqrt2E_{6,{\bf v}'})\},\ 3^7-1+\mbox{$\frac{1}{9}$}\times(82)^2\big\}
\end{eqnarray*}
and $\sqrt2E_7$ is the unique lattice to attain $\kappa^*(B^7)$ (see \cite{Zo01}), up to rotation and reflection. By Remark 2.1, we have
$$\kappa^*_\alpha(B^7)=2954,$$
and it can be attained if and only if the corresponding lattice is $\sqrt2E_7$, up to rotation and reflection.

As a conclusion, Theorem 1 is proved.

\end{proof}

\medskip
\noindent
{\bf Remark 2.2.} In fact, Watson's \cite{Wa01} discussions implies that, for a packing lattice $\Lambda$ of $B^6$, we have

$$\kappa(B^6,\Lambda)\ {=}\left\{
\begin{aligned}
&72,\quad & {\rm if}&\ \Lambda=\sqrt2E_6, \\
&60,\quad & {\rm if}&\ \Lambda=\sqrt2D_6,\\
&56,\quad & {\rm if}&\ \Lambda=\sqrt2D_{5,{\bf v}_6'},\\
&\kappa\leq50,\quad & {\rm if}&\ \Lambda\supseteq\sqrt2D_5\ {\rm and}\ \Lambda\neq\sqrt2E_6, \sqrt2D_6\ {\rm or}\ \sqrt2D_{5,{\bf v}_6'},\\
&\kappa\leq50,\quad & {\rm if}&\ \Lambda\nsupseteq\sqrt2D_5,\\
\end{aligned}
\right .$$
where $\Lambda\supseteq\sqrt2D_5$ means that there exist a $5$-dimensional sub-space $H$ such that $H\cap\Lambda=\sqrt2D_5$ and $\Lambda\nsupseteq\sqrt2D_5$ means that $H\cap\Lambda\neq\sqrt2D_5$ holds for all the $5$-dimensional sub-space $H$. For a packing lattice $\Lambda$ of $B^7$, Watson's \cite{Wa01} discussions implies that

$$\kappa(B^7,\Lambda)\ {=}\left\{
\begin{aligned}
&126,\quad &{\rm if}&\ \Lambda=\sqrt2E_7, \\
&92,\quad & {\rm if}&\ \Lambda=\sqrt2D_{6,{\bf v}_7'}, \\
&84,\quad & {\rm if}&\ \Lambda=\sqrt2E_{6,{\bf v}'}\ {\rm or}\ \sqrt2D_7,\\
&\kappa\leq82,\quad & {\rm if}&\ \Lambda\supseteq\sqrt2E_6\ {\rm and}\ \Lambda\neq\sqrt2E_7, \sqrt2D_{6,{\bf v}_7'}\ {\rm or}\ \sqrt2E_{6,{\bf v}'},\\
&\kappa\leq76,\quad & {\rm if}&\ \Lambda\nsupseteq\sqrt2E_6,\ \Lambda\supseteq\sqrt2D_6 \ {\rm and}\ \Lambda\neq\sqrt2D_7, \sqrt2D_{6,{\bf v}_7'},\\
&\kappa\leq82,\quad & {\rm if}&\ \Lambda\nsupseteq\sqrt2E_6,\ \Lambda\nsupseteq\sqrt2D_6\ {\rm and}\ \Lambda\supseteq\sqrt2D_5,\\
&\kappa\leq82,\quad & {\rm if}&\ \Lambda\nsupseteq\sqrt2D_5.\\
\end{aligned}
\right .$$

\medskip
\noindent
{\bf Remark 2.3.} Let $\Lambda$ be the three-dimensional lattice generated by ${\bf a}_1=(-\frac{2}{3}\sqrt6, \frac{2}{3}\sqrt3,0),$ ${\bf a}_2=(\frac{2}{3}\sqrt6, \frac{2}{3}\sqrt3,0)$ and ${\bf a}_3=(0,\frac{2}{3}\sqrt3, \frac{2}{3}\sqrt6)$. When $\alpha=\frac{4}{3}\sqrt6-2,$ one can verify that \\${\rm card}\{X(\alpha,\Lambda)\}=26.$ Combining with Lemma 2.2, for $\frac{4}{3}\sqrt6-2\leq\alpha<2\sqrt3-2$, we have
$$\kappa_{\alpha}^{*}(B^3)=26.$$

\medskip
\noindent
{\bf Remark 1.} Let $\ell (\Lambda)$ denote the length of the shortest non-zero vector of a given $n$-dimensional lattice $\Lambda$ and let $\gamma (n)$ be a function of $n$. To find a non-zero lattice vector ${\bf v}$ satisfying $\| {\bf v}\|\le \gamma (n) \ell (\Lambda)$ is a fundamental problem in lattice-based cryptography (see \cite{Micc02}). It is easy to see that the generalized kissing numbers provide estimations for the number of such vectors.

\medskip
For almost every $n$-dimensional lattice $\Lambda$ which the length of the shortest non-zero lattice vector is 2, we have
$$\kappa(B^n,\Lambda)=2,$$
see \cite{So01}. However, one can easily imagine that, for a positive number $\alpha$, ${\rm card}\{X(\alpha,\Lambda)\}$ are definitely not the case. Therefore, we end this section by the following problems.

\medskip
\noindent
{\bf Problem 1.} {\it For an $n$-dimensional lattice $\Lambda$ which the length of the shortest non-zero lattice vector is 2 and a positive number $\alpha$, give an algorithm to determine ${\rm card}\{X(\alpha,\Lambda)\}$.}

\medskip
\noindent
{\bf Problem 2.} {\it For an $n$-dimensional lattice $\Lambda$ which the length of the shortest non-zero lattice vector is 2 and a positive number $\alpha$, is there exist an algorithm can determine ${\rm card}\{X(\alpha,\Lambda)\}$ in polynomial time?}

\vspace{1cm}
\noindent
{\LARGE\bf 3. Proof of Theorem 2}

\bigskip
\noindent
For two vectors ${\bf v}$ and ${\bf v}'$ which belongs to an $n$-dimensional lattice $\Lambda$, we say ${\bf v}$ are modulo 2 equivalent to ${\bf v}'$, if and only if
$$\mbox{$\frac{1}{2}$}({\bf v}-{\bf v}')\in\Lambda.$$

Define $X_C(\Lambda)=\{{\bf v}:\ ||{\bf v}||_C=2,\ {\bf v}\in\Lambda\}$, where $\Lambda$ is a packing lattice of $C$. It is easy to see that, the modulo 2 equivalent class of $\Lambda$ which contain ${\bf o}$ cannot contain the vectors of $X_C(\Lambda)$.

Define $C_n=\{(v_1,v_2,\ldots,v_n):\ |v_1|+|v_2|+\ldots+|v_n|\leq1\}.$ Usually, $C_n$ is called an $n$-dimensional cross-polytope.

For a packing lattice $\Lambda$ of $C_n$, suppose there exist two vectors ${\bf v}=(v_1, v_2,\ldots,v_n)$ and ${\bf v}'=(v_1', v_2',\ldots,v_n')$ which belongs to $X_{C_n}(\Lambda)$ with ${\bf v}\neq\pm{\bf v}'$, and ${\bf v}, {\bf v}'$ are modulo 2 equivalent. By $\frac{1}{2}({\bf v}+{\bf v}')\in X_{C_n}(\Lambda)$, we have
$$|v_1+v_1'|+|v_2+v_2'|+\ldots+|v_n+v_n'|=4.$$
Since $|v_1|+|v_2|+\ldots+|v_n|=|v_1'|+|v_2'|+\ldots+|v_n'|=2$, we have
$$v_1v_1',\ v_2v_2',\ \ldots,\ v_nv_n'\geq0.$$
Similarly, by $\frac{1}{2}({\bf v}-{\bf v}')\in X_{C_n}(\Lambda)$, we have
$$v_1v_1',\ v_2v_2',\ \ldots,\ v_nv_n'\leq0.$$

Therefore, we have
$$v_1v_1'=v_2v_2'= \ldots= v_nv_n'=0.\eqno(3.1)$$
In other words, a modulo 2 equivalent class can contain at most $n$ pairs vectors of $X_{C_n}(\Lambda)$.

\medskip
\noindent
{\bf Theorem 2.} {\it Let $C_3$ denote a three-dimensional cross-polytope, then we have
$$\kappa^*(C_3)=18,\quad \kappa^\star(C_3)=14.$$}
\begin{proof}
For a packing lattice $\Lambda$ of $C_3$, if for each modulo 2 equivalent class of $\Lambda$, it contain at most one pair vectors of $X_{C_3}(\Lambda)$, then we have
$${\rm card}\{X_{C_3}(\Lambda)\}\leq(2^3-1)\times2=14,$$
and the equality can be attained. For instance, let $\Lambda$ be the lattice which generated by $\{(-\frac{2}{3}, 1, \frac{1}{3}), (\frac{1}{3}, -\frac{2}{3}, 1), (1, \frac{1}{3}, -\frac{2}{3})\}$, the densest packing lattice of $C_3$ (see \cite{Mi01}), then one can verify that ${\rm card}\{X_{C_3}(\Lambda)\}=14.$

From now on, we suppose there exist a modulo 2 equivalent class of $\Lambda$ which can contain two pairs vectors $\pm{\bf v}_1, \pm{\bf v}_2$ of $X_{C_3}(\Lambda)$. By (3.1), we assume
$${\bf v}_1=(2,0,0),\quad {\bf v}_2=(0,v_{22},2-v_{22}),$$
$$\mbox{$\frac{1}{2}$}({\bf v}_1+{\bf v}_2)=(1,\mbox{$\frac{1}{2}$}v_{22},1-\mbox{$\frac{1}{2}$}v_{22}),\quad \mbox{$\frac{1}{2}$}({\bf v}_1-{\bf v}_2)=(1,-\mbox{$\frac{1}{2}$}v_{22},\mbox{$\frac{1}{2}$}v_{22}-1),$$
where $0\leq v_{22}\leq2$, without loss of generality.

\medskip
\noindent
{\bf Case 1.} There exist a lattice vector ${\bf v}_3\in X_{C_3}(\Lambda)$ with ${\bf v}_3\neq\pm{\bf v}_1, \pm{\bf v}_2$, which are modulo 2 equivalent to ${\bf v}_1$ and ${\bf v}_2$.

By (3.1), we have ${\bf v}_1=(2,0,0),\ {\bf v}_2=(0,2,0),\ {\bf v}_3=(0,0,2)$
and $\pm\frac{1}{2}({\bf v}_{i_1}\pm{\bf v}_{i_2})\in X_{C_3}(\Lambda)$, for $\{i_1, i_2\}\subset \{1,2,3\}$.

If $X_{C_3}(\Lambda)=\{\pm{\bf v}_1, \pm{\bf v}_2, \pm{\bf v}_3, \pm\frac{1}{2}({\bf v}_{i_1}\pm{\bf v}_{i_2})\}$, then we have
$${\rm card}\{X_{C_3}(\Lambda)\}=18.$$
Otherwise, suppose there exist another ${\bf v}=(v_1, v_2, v_3)\in X_{C_3}(\Lambda)$ with
$$v_1, v_2, v_3\geq0,\quad v_1+v_2+v_3=2,\eqno(3.2)$$
without loss of generality. Since ${\bf v}_i-{\bf v}\in\Lambda$, we have
$$v_i\leq 1,\quad i=1,2,3.$$
Combining with $\frac{1}{2}({\bf v}_{i_1}+{\bf v}_{i_2})-{\bf v}\in \Lambda$, we have
$$v_{i_1}+v_{i_2}\leq1$$
holds for $\{i_1, i_2\}\subset \{1,2,3\}$, which contradicts (3.2).

Therefore, in this case, the lattice $\Lambda$ is generated by
$${\bf v}_1=(2,0,0),\quad \mbox{$\frac{1}{2}$}({\bf v}_1+{\bf v}_2)=(1,1,0),\quad \mbox{$\frac{1}{2}$}({\bf v}_1+{\bf v}_3)=(1,0,1)$$
and
$${\rm card}\{X_{C_3}(\Lambda)\}=18.$$

\medskip
\noindent
{\bf Case 2.} For each modulo 2 equivalent classes of $\Lambda$, it contain at most two pairs vectors of $X_{C_3}(\Lambda)$, and there exist another modulo 2 equivalent class of $\Lambda$ which contain two pairs vectors $\pm{\bf v}_1', \pm{\bf v}_2'$ of $X_{C_3}(\Lambda)$.

By (3.1), we may assume
$${\bf v}_1'=(0,2,0),\quad {\bf v}_2'=(v_{21}',0,v_{23}'),$$
$$\mbox{$\frac{1}{2}$}({\bf v}_1'+{\bf v}_2')=(\mbox{$\frac{1}{2}$}v_{21}',1,\mbox{$\frac{1}{2}$}v_{23}'),\quad \mbox{$\frac{1}{2}$}({\bf v}_1'-{\bf v}_2')=(-\mbox{$\frac{1}{2}$}v_{21}',1,-\mbox{$\frac{1}{2}$}v_{23}'),$$
where $|v_{21}'|+|v_{23}'|=2$, without loss of generality.

Since $||{\bf v}_1\pm{\bf v}_2'||_{C_3}\geq2,\ \ ||{\bf v}_1'-{\bf v}_2||_{C_3}\geq2$, we have
$$|v_{21}'|\leq1,\quad v_{22}\leq1.$$
Since $||{\bf v}_2-\frac{1}{2}({\bf v}_1'\pm{\bf v}_2')||_{C_3}\geq2$, by some computation, we have
$$v_{22}+|\mbox{$\frac{1}{2}$}v_{23}'|\leq1.$$
Similarly, since $||{\bf v}_2'\pm\frac{1}{2}({\bf v}_1+{\bf v}_2)||_{C_3}\geq2$, $||{\bf v}_2'\pm\frac{1}{2}({\bf v}_1-{\bf v}_2)||_{C_3}\geq2$, by some computation, we have
$$\mbox{$\frac{1}{2}$}v_{22}\geq|v_{21}'|.$$
Therefore, we have $v_{21}'=0,\ v_{22}=0,\ |v_{23}'|=2$.

In this case, the lattice $\Lambda$ is generated by $\{(2,0,0),(1,1,0),(1,0,1)\}$, which are equivalent to Case 1. In other words, if for each modulo 2 equivalent classes of $\Lambda$, it contain at most two pairs vectors of $X_{C_3}(\Lambda)$, then there are at most one equivalent class which can contain two pairs vectors of $X_{C_3}(\Lambda)$.

\medskip
\noindent
{\bf Case 3.} Besides the modulo 2 equivalent class of $\Lambda$ which contain $\pm{\bf v}_1$ and $\pm{\bf v}_2$, all the modulo 2 equivalent classes of $\Lambda$ contain at most one pair vectors of $X_{C_3}(\Lambda)$.

In this case, we have ${\rm card}\{X_{C_3}(\Lambda)\}\leq16.$ Suppose ${\rm card}\{X_{C_3}(\Lambda)\}=16$, then there exist ${\bf v}_3=(v_{31},v_{32},v_{33})\in X_{C_3}(\Lambda)$ which satisfy:
$${\bf v}_3+\mbox{$\frac{1}{2}$}({\bf v}_1+{\bf v}_2),\quad {\bf v}_3+\mbox{$\frac{1}{2}$}({\bf v}_1-{\bf v}_2),\quad {\bf v}_3+{\bf v}_1\in X_{C_3}(\Lambda).$$
In other words, we have
\begin{eqnarray*}
\hspace{-0.2cm} &| \hspace{-0.35cm}&v_{31}|+|v_{32}|+|v_{33}|=2,\quad |v_{31}+1|+|v_{32}+\mbox{$\frac{1}{2}$}v_{22}|+|v_{33}+1-\mbox{$\frac{1}{2}$}v_{22}|=2,\\
\hspace{-0.2cm} &| \hspace{-0.35cm}&v_{31}+1|+|v_{32}-\mbox{$\frac{1}{2}$}v_{22}|+|v_{33}-(1-\mbox{$\frac{1}{2}$}v_{22})|=2,\quad |v_{31}+2|+|v_{32}|+|v_{33}|=2.
\end{eqnarray*}

By routine computation, we have $v_{31}=-1$ and
$$2|v_{32}|+2-v_{22}=4\quad {\rm or}\quad 2|v_{33}|+v_{22}=4,$$
both leads to a contradiction, which is the lattice $\Lambda$ is generated by $\{(2,0,0),(1,1,0),\\(1,0,1)\}$, which are equivalent to Case 1. Therefore, in this case, we have
$${\rm card}\{X_{C_3}(\Lambda)\}\leq14.$$

\medskip
As a conclusion of all the cases, we have
$$\kappa^*(C_3)=18,\quad \kappa^\star(C_3)=14.$$
Theorem 2 is proved.
\end{proof}

\vspace{1cm}
\noindent
{\LARGE\bf 4. Proof of Theorem 3}

\bigskip
\noindent
{\bf Lemma 4.1.} {\it If there exist two vectors ${\bf v}_1, {\bf v}_2\in X_{C_4}(\Lambda)$ are modulo 2 equivalent , two vectors ${\bf v}_3, {\bf v}_4\in X_{C_4}(\Lambda)$ are modulo 2 equivalent, ${\bf v}_i\neq\pm{\bf v}_j$, $\{i,j\}\subset\{1,2,3,4\}$ and
$$\mbox{$\frac{1}{2}$}({\bf v}_1+{\bf v}_2)=\mbox{$\frac{1}{2}$}({\bf v}_3+{\bf v}_4),$$
denote it by ${\bf v}=(v_1, v_2, v_3, v_4)$. Then we have $|v_1|=|v_2|=|v_3|=|v_4|=\frac{1}{2}$.}

\begin{proof}
Denote ${\bf v}_1=(v_{11}, v_{12}, v_{13}, v_{14})$, ${\bf v}_2=(v_{21}, v_{22}, v_{23}, v_{24})$ with
$$|v_{11}|+|v_{12}|+|v_{13}|+|v_{14}|=|v_{21}|+|v_{22}|+|v_{23}|+|v_{24}|=2.$$

If one of $v_{11}, v_{12}, v_{13}, v_{14}=\pm2$, for instance, $v_{11}=2$. Combining with (3.1), we have
$${\bf v}_1=(2,0,0,0),\quad {\bf v}_2=(0, v_{22}, v_{23}, v_{24}),\quad \mbox{$\frac{1}{2}$}({\bf v}_1+{\bf v}_2)=(1, \mbox{$\frac{1}{2}$}v_{22}, \mbox{$\frac{1}{2}$}v_{23}, \mbox{$\frac{1}{2}$}v_{24}).$$
To satisfy $\mbox{$\frac{1}{2}$}({\bf v}_1+{\bf v}_2)=\mbox{$\frac{1}{2}$}({\bf v}_3+{\bf v}_4),$ by (3.1), we have one of $\pm{\bf v}_3, \pm{\bf v}_4$ is $\pm(2,0,0,0)$, which is a contradiction.

Therefore, by (3.1), we may assume ${\bf v}_1=(2v_1, 2v_2, 0, 0)$, ${\bf v}_2=(0, 0, 2v_3, 2v_4)$ with $|v_1|+|v_2|=|v_3|+|v_4|=1$ and $\mbox{$\frac{1}{2}$}({\bf v}_1+{\bf v}_2)=(v_1, v_2, v_3, v_4)$, without loss of generality. To satisfy $\mbox{$\frac{1}{2}$}({\bf v}_1+{\bf v}_2)=\mbox{$\frac{1}{2}$}({\bf v}_3+{\bf v}_4),$ by (3.1), we may assume ${\bf v}_3=(2v_1, 0, 2v_3, 0)$, ${\bf v}_4=(0, 2v_2, 0, 2v_4)$ with $|v_1|+|v_3|=|v_2|+|v_4|=1$. Therefore, we have
$$|v_1|=|v_2|=|v_3|=|v_4|=\mbox{$\frac{1}{2}$},$$
Lemma 4.1 is proved.
\end{proof}

\medskip
\noindent
{\bf Remark 4.1.} Suppose the lattice $\Lambda$ satisfies the conditions in Lemma 4.1, without loss of generality we assume $v_1=v_2=v_3=v_4=\frac{1}{2}.$ Then, by (3.1), there exist modulo 2 equivalent classes $U_i(\Lambda)$, $i\in\{1,\ldots,6\}$, which satisfy:
\begin{eqnarray*}
U_1(\Lambda)\cap X_{C_4}(\Lambda)\hspace{-0.2cm} &= \hspace{-0.2cm}&\{\pm{\bf v}_1=\pm(1,1,0,0),\ \pm{\bf v}_2=\pm(0,0,1,1)\},\\
U_2(\Lambda)\cap X_{C_4}(\Lambda)\hspace{-0.2cm} &= \hspace{-0.2cm}&\{\pm{\bf v}_3=\pm(1,0,1,0),\ \pm{\bf v}_4=\pm(0,1,0,1)\},\\
U_3(\Lambda)\cap X_{C_4}(\Lambda)\hspace{-0.2cm} &= \hspace{-0.2cm}&\{\pm({\bf v}_1-{\bf v}_4)=\pm(1,0,0,-1),\ \pm({\bf v}_1-{\bf v}_3)=\pm(0,1,-1,0)\},\\
U_4(\Lambda)\cap X_{C_4}(\Lambda)\hspace{-0.2cm} &= \hspace{-0.2cm}&\{\pm\mbox{$\frac{1}{2}$}({\bf v}_1+{\bf v}_2)=\pm(\mbox{$\frac{1}{2},\frac{1}{2},\frac{1}{2},\frac{1}{2})$}\},\\
U_5(\Lambda)\cap X_{C_4}(\Lambda)\hspace{-0.2cm} &= \hspace{-0.2cm}&\{\pm\mbox{$\frac{1}{2}$}({\bf v}_1-{\bf v}_2)=\pm(\mbox{$\frac{1}{2},\frac{1}{2},-\frac{1}{2},-\frac{1}{2}$})\},\\
U_6(\Lambda)\cap X_{C_4}(\Lambda)\hspace{-0.2cm} &= \hspace{-0.2cm}&\{\pm\mbox{$\frac{1}{2}$}({\bf v}_3-{\bf v}_4)=\pm(\mbox{$\frac{1}{2},-\frac{1}{2},\frac{1}{2},-\frac{1}{2}$})\}.
\end{eqnarray*}

Furthermore, if we have one of $\pm(2,0,0,0), \pm(0,2,0,0), \pm(0,0,2,0), \pm(0,0,0,2)$ belongs to $X_{C_4}(\Lambda)$, then it is easy to verify that, the lattice $\Lambda$ is generated by $\{(2,0,0,0), (1,1,0,0),$ $(1,0,1,0), (\frac{1}{2}, \frac{1}{2}, \frac{1}{2}, \frac{1}{2})\}$ and ${\rm card}\{X_{C_4}(\Lambda)\}=40.$

\medskip
\noindent
{\bf Theorem 3.} {\it Let $C_4$ denote a four-dimensional cross-polytope, then we have
$$\kappa^*(C_4)=40.$$}
\begin{proof}

For a packing lattice $\Lambda$ of $C_4$, if there exist a modulo 2 equivalent class of $\Lambda$ which can contain four pairs vectors $\pm{\bf v}_1, \pm{\bf v}_2, \pm{\bf v}_3, \pm{\bf v}_4$ of $X_{C_4}(\Lambda)$, by (3.1), we set
$${\bf v}_1=(2,0,0,0),\ {\bf v}_2=(0,2,0,0),\ {\bf v}_3=(0,0,2,0),\ {\bf v}_4=(0,0,0,2)$$
and $\pm\frac{1}{2}({\bf v}_{i_1}\pm{\bf v}_{i_2})\in X_{C_4}(\Lambda)$, for $\{i_1, i_2\}\subset \{1,2,3,4\}$.

Suppose $X_{C_4}(\Lambda)=\{\pm{\bf v}_1, \pm{\bf v}_2, \pm{\bf v}_3, \pm{\bf v}_4, \pm\frac{1}{2}({\bf v}_{i_1}\pm{\bf v}_{i_2})\}$, then we have
$${\rm card}\{X_{C_4}(\Lambda)\}=32.$$
Otherwise, suppose there exist another vector ${\bf v}=(v_1, v_2, v_3, v_4)\in X_{C_4}(\Lambda)$ with
$$v_1, v_2, v_3, v_4\geq0,\quad v_1+v_2+v_3+v_4=2,\eqno(4.1)$$
without loss of generality. Since ${\bf v}_i-{\bf v}\in\Lambda$, we have
$$v_i\leq 1,\quad i=1,2,3,4.$$
Combining with $\frac{1}{2}({\bf v}_{i_1}+{\bf v}_{i_2})-{\bf v}\in \Lambda$, we have
$$v_{i_1}+v_{i_2}\leq1$$
holds for $\{i_1, i_2\}\subset \{1,2,3,4\}$. Combining with (4.1), we have
$$v_1=v_2=v_3=v_4=\mbox{$\frac{1}{2}$}.$$
In this case, the lattice $\Lambda$ is generated by ${\bf v}_1=(2,0,0,0),\ \mbox{$\frac{1}{2}$}({\bf v}_1+{\bf v}_2)=(1,1,0,0),\ \mbox{$\frac{1}{2}$}({\bf v}_1+{\bf v}_3)=(1,0,1,0),\ \mbox{$\frac{1}{4}$}({\bf v}_1+{\bf v}_2+{\bf v}_3+{\bf v}_4)=(\mbox{$\frac{1}{2},\frac{1}{2},\frac{1}{2},\frac{1}{2}$})$
and
$${\rm card}\{X_{C_4}(\Lambda)\}=40.$$
Since $$\left |\begin{array}{cccc}\medskip
2 & 0 & 0 & 0\\
\medskip
1 & 1 & 0 & 0\\
\medskip
1 & 0 & 1 & 0\\
\frac{1}{2} & \frac{1}{2} & \frac{1}{2} & \frac{1}{2}\\
\end{array}\right|=1,$$
combining with ${\rm vol}(C_4)=\frac{2^4}{4!}=\frac{2}{3}$, we have the packing density of $C_4+\Lambda$ is $2/3$.

If there exist two modulo 2 equivalent class which contain exactly three pairs vectors $\pm{\bf v}_1,\ \pm{\bf v}_2,\ \pm{\bf v}_3$ and $\pm{\bf v}_1',\ \pm{\bf v}_2',\ \pm{\bf v}_3'$ of $X_{C_4}(\Lambda)$, by (3.1), we set
$${\bf v}_1=(2,0,0,0),\quad {\bf v}_2=(0,2,0,0),\quad {\bf v}_3=(0,0,v_{33},v_{34}),$$
$${\bf v}_1'=(0,0,2,0),\quad {\bf v}_2'=(0,0,0,2),\quad {\bf v}_3'=(v_{31}',v_{32}',0,0),$$
where $v_{31}',\ v_{32}',\ v_{33},\ v_{34}\geq0$, $v_{31}'+v_{32}'=2,\ v_{33}+v_{34}=2$, without loss of generality.

Since ${\bf v}_1-{\bf v}_3',\ {\bf v}_2-{\bf v}_3'\in\Lambda$, we have $v_{31}'=v_{32}'=1$. Similarly, by ${\bf v}_1'-{\bf v}_3, {\bf v}_2'-{\bf v}_3\in\Lambda$, we have $v_{33}=v_{34}=1$.

Since $\frac{1}{2}({\bf v}_1+{\bf v}_3), \frac{1}{2}({\bf v}'_1+{\bf v}'_3)\in\Lambda$, we have
$$\mbox{$\frac{1}{2}({\bf v}_1'+{\bf v}_3')-\frac{1}{2}({\bf v}_1+{\bf v}_3)=(-\frac{1}{2}, \frac{1}{2}, \frac{1}{2}, -\frac{1}{2}$})\in\Lambda.$$
Combining with $\frac{1}{2}({\bf v}_1'-{\bf v}_2')\in\Lambda$, we have
$$\mbox{$\frac{1}{2}({\bf v}_1'+{\bf v}_3')-\frac{1}{2}({\bf v}_1+{\bf v}_3)-\frac{1}{2}({\bf v}_1'-{\bf v}_2')=(-\frac{1}{2}, \frac{1}{2}, -\frac{1}{2}, \frac{1}{2}$})\in\Lambda.$$
Since $$\left |\begin{array}{cccc}\medskip
-\frac{1}{2} & -\frac{1}{2} & 1 & 0\\
\medskip
\frac{1}{2} & \frac{1}{2} & 1 & 0\\
\medskip
-\frac{1}{2} & \frac{1}{2} & \frac{1}{2} & -\frac{1}{2}\\
-\frac{1}{2} & \frac{1}{2} & -\frac{1}{2} & \frac{1}{2}\\
\end{array}\right|=1,$$
where $(-\frac{1}{2}, -\frac{1}{2}, 1, 0)=\frac{1}{2}({\bf v}_1'-{\bf v}_3')\in\Lambda$, $(\frac{1}{2}, \frac{1}{2}, 1, 0)=\frac{1}{2}({\bf v}_1'+{\bf v}_3')\in\Lambda$, we have
$$\mbox{$(-\frac{1}{2}, -\frac{1}{2}, 1, 0),\quad (\frac{1}{2}, \frac{1}{2}, 1, 0),\quad (-\frac{1}{2}, \frac{1}{2}, \frac{1}{2}, -\frac{1}{2}),\quad (-\frac{1}{2}, \frac{1}{2}, -\frac{1}{2}, \frac{1}{2})$}$$
is a basis of $\Lambda$ and the packing density of $C_4+\Lambda$ is $2/3$. By enumerate, in this case, we have
$${\rm card}\{X_{C_4}(\Lambda)\}=36.$$

On the other hand, if for each modulo 2 equivalent classes of $\Lambda$, it can contain at most one pair vectors of $X_{C_4}(\Lambda)$, then we have
$${\rm card}\{X_{C_4}(\Lambda)\}\leq2\times(2^4-1)=30.$$

\medskip
From now on, we deal with the packing lattice $\Lambda$ of $C_4$ such that, there exist a modulo 2 equivalent class of $\Lambda$ which can contain at least two pairs vectors of $X_{C_4}(\Lambda)$, but for each modulo 2 equivalent classes of $\Lambda$, it can contain at most three pairs vectors of $X_{C_4}(\Lambda)$ and there are at most one modulo 2 equivalent class of $\Lambda$ can contain three pairs vectors of $X_{C_4}(\Lambda)$.

\medskip
\noindent
{\bf Case 1.} There exist two vectors ${\bf v}_1, {\bf v}_2\in X_{C_4}(\Lambda)$ are modulo 2 equivalent, two vectors ${\bf v}_3, {\bf v}_4\in X_{C_4}(\Lambda)$ are modulo 2 equivalent, ${\bf v}_i\neq\pm{\bf v}_j$, $\{i,j\}\subset\{1,2,3,4\}$ and $\frac{1}{2}({\bf v}_1+{\bf v}_2)=\frac{1}{2}({\bf v}_3+{\bf v}_4)$.

In this case, we use the notations and results of Remark 4.1.

\medskip
\noindent
{\bf Case 1.1.} One of $\pm(2,0,0,0), \pm(0,2,0,0), \pm(0,0,2,0), \pm(0,0,0,2)$ belongs to $\Lambda$.

By Remark 4.1, we have
$${\rm card}\{X_{C_4}(\Lambda)\}=40.$$

\medskip
\noindent
{\bf Case 1.2.} Besides $U_i(\Lambda)$, $i\in\{1,\ldots,6\}$ of Remark 4.1, all the modulo 2 equivalent classes of $\Lambda$ contain at most one pair vectors of $X_{C_4}(\Lambda)$.

Combining with Remark 4.1, we have
$${\rm card}\{X_{C_4}(\Lambda)\}\leq18+(15-6)\times2=36.$$

\medskip
\noindent
{\bf Case 1.3.} Besides $U_i(\Lambda)$, $i\in\{1,\ldots,6\}$ of Remark 4.1, there exist another modulo 2 equivalent class of $\Lambda$, which can contain two pairs vectors $\pm{\bf v}, \pm{\bf v}'$ of $X_{C_4}(\Lambda)$ and $\pm(2,0,0,0), \pm(0,2,0,0), \pm(0,0,2,0), \pm(0,0,0,2)\notin\Lambda$.

By (3.1), we set
$${\bf v}=(v_1,v_2,0,0),\quad {\bf v}'=(0,0,v_3',v_4'),$$
where $|v_1|+|v_2|=2$, $|v_3'|+|v_4'|=2$ and $v_1,\ v_3'>0$, without loss of generality.

By ${\bf v}_1-{\bf v}\in\Lambda,\ {\bf v}_2-{\bf v}'\in\Lambda$, we have $v_2<0,\ v_4'<0$, which means $v_1-v_2=2,\ v_3'-v_4'=2$. Therefore, we have
$$|1-v_1|+|1+v_2|=v_1+v_2\quad {\rm or}\quad -(v_1+v_2),\eqno(4.2)$$
$$|1-v_3'|+|1+v_4'|=v_3'+v_4'\quad {\rm or}\quad -(v_3'+v_4').\eqno(4.3)$$
However, since $\frac{1}{2}({\bf v}+{\bf v}')=(\mbox{$\frac{1}{2}$}v_{1},\mbox{$\frac{1}{2}$}v_{2},\mbox{$\frac{1}{2}$}v_{3}',\mbox{$\frac{1}{2}$}v_{4}')\in\Lambda$ and $\frac{1}{2}({\bf v}_3-{\bf v}_4)=(\frac{1}{2}, -\frac{1}{2}, \frac{1}{2}, -\frac{1}{2})\in\Lambda$, we have
$$|\mbox{$\frac{1}{2}$}-\mbox{$\frac{1}{2}$}v_{1}|+|\mbox{$\frac{1}{2}$}+\mbox{$\frac{1}{2}$}v_{2}|+|\mbox{$\frac{1}{2}$}-\mbox{$\frac{1}{2}$}v_{3}'|+|\mbox{$\frac{1}{2}$}+\mbox{$\frac{1}{2}$}v_{4}'|\geq2.$$
Combining with (4.2) and (4.3), we have $\{|v_1|, |v_2|\}=\{|v_3'|, |v_4'|\}=\{0,2\}$, namely:
$${\bf v},\ {\bf v}'\in \{\pm(2,0,0,0),\ \pm(0,2,0,0),\ \pm(0,0,2,0),\ \pm(0,0,0,2)\},$$
which is a contradiction.

\medskip
From now on, we further assume that, for arbitrary two modulo 2 equivalent vectors ${\bf v}, {\bf v}'\in X_{C_4}(\Lambda)$ with ${\bf v}\neq\pm{\bf v}'$, the mid-point of ${\bf v}$ and ${\bf v}'$ cannot be the mid-point of other two different vectors of $X_{C_4}(\Lambda)$.

\medskip
\noindent
{\bf Case 2.} For arbitrary two modulo 2 equivalent vectors ${\bf v}, {\bf v}'\in X_{C_4}(\Lambda)$ with ${\bf v}\neq\pm{\bf v}'$, the mid-point of ${\bf v}$ and ${\bf v}'$ cannot modulo 2 equivalent to vectors (besides $\pm\frac{1}{2}({\bf v}+{\bf v}')$) of $X_{C_4}(\Lambda)$.

\medskip
\noindent
{\bf Case 2.1.} All the modulo 2 equivalent classes of $\Lambda$ contain at most two pairs vectors of $X_{C_4}(\Lambda)$.

\medskip
For two modulo 2 equivalent vectors ${\bf v}, {\bf v}'\in X_{C_4}(\Lambda)$ with ${\bf v}\neq\pm{\bf v}'$, since the mid-point of ${\bf v}$ and ${\bf v}'$ and the mid-point of ${\bf v}$ and $(-{\bf v}')$ cannot equivalent to other vectors of $X_{C_4}(\Lambda)$, there exist modulo 2 equivalent classes $U_i(\Lambda)$, $i\in\{1,2,3\}$, which satisfy:
\begin{eqnarray*}
U_1(\Lambda)\cap X_{C_4}(\Lambda)\hspace{-0.2cm} &= \hspace{-0.2cm}&\{\pm{\bf v}, \pm{\bf v}'\},\quad U_2(\Lambda)\cap X_{C_4}(\Lambda)=\{\pm\mbox{$\frac{1}{2}$}({\bf v}+{\bf v}')\},\\
U_3(\Lambda)\cap X_{C_4}(\Lambda)\hspace{-0.2cm} &= \hspace{-0.2cm}&\{\pm\mbox{$\frac{1}{2}$}({\bf v}-{\bf v}')\}.
\end{eqnarray*}
Combining with $\frac{1}{2}({\bf v}+{\bf v}')$ and $\frac{1}{2}({\bf v}-{\bf v}')$ cannot be the mid-point of other two vectors of $X_{C_4}(\Lambda)$, by induction we have
$${\rm card}\{X_{C_4}(\Lambda)\}\leq(4+2+2)\times5=40.$$

\medskip
\noindent
{\bf Case 2.2.} There exist a modulo 2 equivalent class of $\Lambda$, which can contain three pairs vectors $\pm{\bf v}_1, \pm{\bf v}_2, \pm{\bf v}_3$ of $X_{C_4}(\Lambda)$.

In this case, there exist modulo 2 equivalent classes $U_i(\Lambda)$, $i\in\{1,\ldots,7\}$, which satisfy
\begin{eqnarray*}
U_1(\Lambda)\cap X_{C_4}(\Lambda)\hspace{-0.2cm} &= \hspace{-0.2cm}&\{\pm{\bf v}_1,\ \pm{\bf v}_2,\ \pm{\bf v}_3\},\\
U_2(\Lambda)\cap X_{C_4}(\Lambda)\hspace{-0.2cm} &= \hspace{-0.2cm}&\{\pm\mbox{$\frac{1}{2}$}({\bf v}_1+{\bf v}_2)\},\quad U_3(\Lambda)\cap X_{C_4}(\Lambda)=\{\pm\mbox{$\frac{1}{2}$}({\bf v}_1-{\bf v}_2)\},\\
U_4(\Lambda)\cap X_{C_4}(\Lambda)\hspace{-0.2cm} &= \hspace{-0.2cm}&\{\pm\mbox{$\frac{1}{2}$}({\bf v}_1+{\bf v}_3)\},\quad U_5(\Lambda)\cap X_{C_4}(\Lambda)=\{\pm\mbox{$\frac{1}{2}$}({\bf v}_1-{\bf v}_3)\},\\
U_6(\Lambda)\cap X_{C_4}(\Lambda)\hspace{-0.2cm} &= \hspace{-0.2cm}&\{\pm\mbox{$\frac{1}{2}$}({\bf v}_2+{\bf v}_3)\},\quad U_7(\Lambda)\cap X_{C_4}(\Lambda)=\{\pm\mbox{$\frac{1}{2}$}({\bf v}_2-{\bf v}_3)\}.
\end{eqnarray*}
Since all the rest of modulo 2 equivalent classes of $\Lambda$ (besides $U_1(\Lambda)$) can contain at most two pairs vectors of $X_{C_4}(\Lambda)$, we have
$${\rm card}\{X_{C_4}(\Lambda)\}\leq18+(4+2+2)\times2+2+2=38.$$

\medskip
\noindent
{\bf Case 3.} There exist two modulo 2 equivalent vectors ${\bf v}_1,\ {\bf v}_2\in X_{C_4}(\Lambda)$ with ${\bf v}_1\neq\pm{\bf v}_2$ and there exist ${\bf v}_3\in X_{C_4}(\Lambda)$ are modulo 2 equivalent to $\frac{1}{2}({\bf v}_1+{\bf v}_2)$ with ${\bf v}_3\neq\pm\frac{1}{2}({\bf v}_1+{\bf v}_2).$

\medskip
\noindent
{\bf Case 3.1.} All the modulo 2 equivalent classes of $\Lambda$ contain at most two pairs vectors of $X_{C_4}(\Lambda)$.

\medskip
By (3.1), we set
$${\bf v}_1=(0,2,0,0),\quad {\bf v}_2=(0,0,v_{23},v_{24}),$$
where $v_{23},\ v_{24}\geq0$, $v_{23}+v_{24}=2$, without loss of generality. We have one of the following scenarios must happen.

\medskip
\noindent
{\bf Scenario (A).} Both the mid-point of ${\bf v}_1$ and $(-{\bf v}_2)$, the mid-point of $\frac{1}{2}({\bf v}_1+{\bf v}_2)$ and ${\bf v}_3$ and the mid-point of $\frac{1}{2}({\bf v}_1+{\bf v}_2)$ and $(-{\bf v}_3)$ are not equivalent to vectors (besides $\pm\frac{1}{2}({\bf v}_1-{\bf v}_2)$, $\pm\frac{1}{2}(\frac{1}{2}({\bf v}_1+{\bf v}_2)+{\bf v}_3)$ and $\pm\frac{1}{2}(\frac{1}{2}({\bf v}_1+{\bf v}_2)-{\bf v}_3)$, respectively) of $X_{C_4}(\Lambda)$.

Suppose $v_{23}\cdot v_{24}\neq0$, by (3.1) we have $\pm{\bf v}_3=\pm(2,0,0,0)$ are modulo 2 equivalent to $\pm\frac{1}{2}({\bf v}_1+{\bf v}_2)=\pm(0,1,\mbox{$\frac{1}{2}$}v_{23},\mbox{$\frac{1}{2}$}v_{24})$. Therefore, there exist modulo 2 equivalent classes $U_i(\Lambda)$, $i\in\{1,\ldots,5\}$, which satisfy:
\begin{eqnarray*}
U_1(\Lambda)\cap X_{C_4}(\Lambda)\hspace{-0.2cm} &= \hspace{-0.2cm}&\{\pm{\bf v}_1=\pm(0,2,0,0),\ \pm{\bf v}_2=\pm(0,0,v_{23},v_{24})\},\\
U_2(\Lambda)\cap X_{C_4}(\Lambda)\hspace{-0.2cm} &= \hspace{-0.2cm}&\{\pm\mbox{$\frac{1}{2}$}({\bf v}_1+{\bf v}_2)=\pm(0,1,\mbox{$\frac{1}{2}$}v_{23},\mbox{$\frac{1}{2}$}v_{24}),\ \pm{\bf v}_3=\pm(2,0,0,0)\},\\
U_3(\Lambda)\cap X_{C_4}(\Lambda)\hspace{-0.2cm} &= \hspace{-0.2cm}&\{\pm\mbox{$\frac{1}{2}$}({\bf v}_1-{\bf v}_2)=\pm(0,1,-\mbox{$\frac{1}{2}$}v_{23},-\mbox{$\frac{1}{2}$}v_{24})\},\\
U_4(\Lambda)\cap X_{C_4}(\Lambda)\hspace{-0.2cm} &= \hspace{-0.2cm}&\{\pm\mbox{$\frac{1}{2}(\frac{1}{2}({\bf v}_1+{\bf v}_2)+{\bf v}_3)=\pm(1,\frac{1}{2},\mbox{$\frac{1}{4}$}v_{23},\mbox{$\frac{1}{4}$}v_{24}$})\},\\
U_5(\Lambda)\cap X_{C_4}(\Lambda)\hspace{-0.2cm} &= \hspace{-0.2cm}&\{\pm\mbox{$\frac{1}{2}(\frac{1}{2}({\bf v}_1+{\bf v}_2)-{\bf v}_3)=\pm(-1,\frac{1}{2},\mbox{$\frac{1}{4}$}v_{23},\mbox{$\frac{1}{4}$}v_{24}$})\},
\end{eqnarray*}
up to permutation of coordinates and signs of coordinates, denote it by scenario (${\rm A_1}$).

On the contrary, suppose $v_{23}=0, v_{24}=2$. By (3.1), we set $\pm{\bf v}_3=\pm(v_{31},0,v_{33},0)$ are modulo 2 equivalent to $\pm\frac{1}{2}({\bf v}_1+{\bf v}_2)=\pm(0,1,0,1)$. Therefore, there exist modulo 2 equivalent classes $U_i(\Lambda)$, $i\in\{1,\ldots,5\}$, which satisfy:
\begin{eqnarray*}
U_1(\Lambda)\cap X_{C_4}(\Lambda)\hspace{-0.2cm} &= \hspace{-0.2cm}&\{\pm{\bf v}_1=\pm(0,2,0,0),\ \pm{\bf v}_2=\pm(0,0,0,2)\},\\
U_2(\Lambda)\cap X_{C_4}(\Lambda)\hspace{-0.2cm} &= \hspace{-0.2cm}&\{\pm\mbox{$\frac{1}{2}$}({\bf v}_1+{\bf v}_2)=\pm(0,1,0,1),\ \pm{\bf v}_3=\pm(v_{31},0,v_{33},0)\},\\
U_3(\Lambda)\cap X_{C_4}(\Lambda)\hspace{-0.2cm} &= \hspace{-0.2cm}&\{\pm\mbox{$\frac{1}{2}$}({\bf v}_1-{\bf v}_2)=\pm(0,1,0,-1)\},\\
U_4(\Lambda)\cap X_{C_4}(\Lambda)\hspace{-0.2cm} &= \hspace{-0.2cm}&\{\pm\mbox{$\frac{1}{2}(\frac{1}{2}({\bf v}_1+{\bf v}_2)+{\bf v}_3)=\pm(\mbox{$\frac{1}{2}$}v_{31},\frac{1}{2},\mbox{$\frac{1}{2}$}v_{33},\frac{1}{2}$})\},\\
U_5(\Lambda)\cap X_{C_4}(\Lambda)\hspace{-0.2cm} &= \hspace{-0.2cm}&\{\pm\mbox{$\frac{1}{2}(\frac{1}{2}({\bf v}_1+{\bf v}_2)-{\bf v}_3)=\pm(-\mbox{$\frac{1}{2}$}v_{31},\frac{1}{2},-\mbox{$\frac{1}{2}$}v_{33},\frac{1}{2}$})\},
\end{eqnarray*}
up to permutation of coordinates and signs of coordinates, denote it by scenario (${\rm A_2}$).

\medskip
\noindent
{\bf Scenario (B).} There exist $\pm{\bf v}_4\in X_{C_4}(\Lambda)$, which are modulo 2 equivalent to $\pm\frac{1}{2}({\bf v}_1-{\bf v}_2)$. (${\bf v}_4\neq\pm\frac{1}{2}({\bf v}_1-{\bf v}_2)$)

In this scenario, if $v_{23}\cdot v_{24}\neq0$, by (3.1) we have both $\pm{\bf v}_3=\pm(2,0,0,0)$ and $\pm{\bf v}_4=\pm(2,0,0,0)$, which is a contradiction. Therefore, we set $v_{23}=0, v_{24}=2$.

By (3.1), we assume $\pm{\bf v}_3=\pm(v_{31},0,v_{33},0)$ are modulo 2 equivalent to $\pm\mbox{$\frac{1}{2}$}({\bf v}_1+{\bf v}_2)$ $=\pm(0,1,0,1)$ and $\pm{\bf v}_4=\pm(v_{41},0,v_{43},0)$ are modulo 2 equivalent to $\pm\mbox{$\frac{1}{2}$}({\bf v}_1-{\bf v}_2)=\pm(0,1,0,-1)$ with $|v_{31}|+|v_{33}|=2,\ |v_{41}|+|v_{43}|=2.$

Suppose there exist ${\bf v}_5\in X_{C_4}(\Lambda)$ (${\bf v}_5\neq\pm\frac{1}{2}(\frac{1}{2}({\bf v}_1+{\bf v}_2)+{\bf v}_3)$) are modulo 2 equivalent to $\pm\frac{1}{2}(\frac{1}{2}({\bf v}_1+{\bf v}_2)+{\bf v}_3)=\pm(\mbox{$\frac{1}{2}$}v_{31},\frac{1}{2},\mbox{$\frac{1}{2}$}v_{33},\frac{1}{2})$, by (3.1) we assume $v_{31}=2$, $v_{33}=0$ and ${\bf v}_5=(0,0,2,0)$, without loss of generality.

Since $||{\bf v}_3\pm{\bf v}_4||_{C_4},\ ||{\bf v}_5\pm{\bf v}_4||_{C_4}\geq2$, we have $|v_{41}|=|v_{43}|=1$. However, since
\begin{eqnarray*}
\pm\mbox{$\frac{1}{2}(\frac{1}{2}(\frac{1}{2}$}({\bf v}_1+{\bf v}_2)+{\bf v}_3)+{\bf v}_5)\hspace{-0.2cm} &= \hspace{-0.2cm}&\pm(\mbox{$\frac{1}{2},\frac{1}{4},1,\frac{1}{4}$})\in\Lambda,\\
\pm\mbox{$\frac{1}{2}(\frac{1}{2}(\frac{1}{2}$}({\bf v}_1+{\bf v}_2)+{\bf v}_3)-{\bf v}_5)\hspace{-0.2cm} &= \hspace{-0.2cm}&\pm(\mbox{$\frac{1}{2},\frac{1}{4},-1,\frac{1}{4}$})\in\Lambda,
\end{eqnarray*}
we have
$$\min\big\{||\mbox{$\frac{1}{2}(\frac{1}{2}(\frac{1}{2}$}({\bf v}_1+{\bf v}_2)+{\bf v}_3)+{\bf v}_5)\pm{\bf v}_4||_{C_4},\ ||\mbox{$\frac{1}{2}(\frac{1}{2}(\frac{1}{2}$}({\bf v}_1+{\bf v}_2)+{\bf v}_3)-{\bf v}_5)\pm{\bf v}_4||_{C_4}\big\}=1,$$
which is a contradiction.

Therefore, there exist modulo 2 equivalent classes $U_i(\Lambda)$, $i\in\{1,\ldots,7\}$, which satisfy:
\begin{eqnarray*}
U_1(\Lambda)\cap X_{C_4}(\Lambda)\hspace{-0.2cm} &= \hspace{-0.2cm}&\{\pm{\bf v}_1=\pm(0,2,0,0),\ \pm{\bf v}_2=\pm(0,0,0,2)\},\\
U_2(\Lambda)\cap X_{C_4}(\Lambda)\hspace{-0.2cm} &= \hspace{-0.2cm}&\{\pm\mbox{$\frac{1}{2}$}({\bf v}_1+{\bf v}_2)=\pm(0,1,0,1),\ \pm{\bf v}_3=\pm(v_{31},0,v_{33},0)\},\\
U_3(\Lambda)\cap X_{C_4}(\Lambda)\hspace{-0.2cm} &= \hspace{-0.2cm}&\{\pm\mbox{$\frac{1}{2}$}({\bf v}_1-{\bf v}_2)=\pm(0,1,0,-1),\ \pm{\bf v}_4=\pm(v_{41},0,v_{43},0)\},\\
U_4(\Lambda)\cap X_{C_4}(\Lambda)\hspace{-0.2cm} &= \hspace{-0.2cm}&\{\pm\mbox{$\frac{1}{2}(\frac{1}{2}({\bf v}_1+{\bf v}_2)+{\bf v}_3)=\pm(\mbox{$\frac{1}{2}$}v_{31},\frac{1}{2},\mbox{$\frac{1}{2}$}v_{33},\frac{1}{2}$})\},\\
U_5(\Lambda)\cap X_{C_4}(\Lambda)\hspace{-0.2cm} &= \hspace{-0.2cm}&\{\pm\mbox{$\frac{1}{2}(\frac{1}{2}({\bf v}_1+{\bf v}_2)-{\bf v}_3)=\pm(-\mbox{$\frac{1}{2}$}v_{31},\frac{1}{2},-\mbox{$\frac{1}{2}$}v_{33},\frac{1}{2}$})\},\\
U_6(\Lambda)\cap X_{C_4}(\Lambda)\hspace{-0.2cm} &= \hspace{-0.2cm}&\{\pm\mbox{$\frac{1}{2}(\frac{1}{2}({\bf v}_1-{\bf v}_2)+{\bf v}_4)=\pm(\mbox{$\frac{1}{2}$}v_{41},\frac{1}{2},\mbox{$\frac{1}{2}$}v_{43},-\frac{1}{2}$})\},\\
U_7(\Lambda)\cap X_{C_4}(\Lambda)\hspace{-0.2cm} &= \hspace{-0.2cm}&\{\pm\mbox{$\frac{1}{2}(\frac{1}{2}({\bf v}_1-{\bf v}_2)-{\bf v}_4)=\pm(-\mbox{$\frac{1}{2}$}v_{41},\frac{1}{2},-\mbox{$\frac{1}{2}$}v_{43},-\frac{1}{2}$})\},
\end{eqnarray*}
up to permutation of coordinates and signs of coordinates.

\medskip
\noindent
{\bf Scenario (C).} The mid-point of ${\bf v}_1$ and $(-{\bf v}_2)$ are not equivalent to vectors (besides $\pm\frac{1}{2}({\bf v}_1-{\bf v}_2)$) of $X_{C_4}(\Lambda)$, and there exist ${\bf v}_4\in X_{C_4}(\Lambda)$ which are modulo 2 equivalent to $\frac{1}{2}(\frac{1}{2}({\bf v}_1+{\bf v}_2)+{\bf v}_3)$. (${\bf v}_4\neq\pm\frac{1}{2}(\frac{1}{2}({\bf v}_1+{\bf v}_2)+{\bf v}_3)$)

Combining with the discussions in scenarios (A) and (B) and (3.1), in this scenario, there exist modulo 2 equivalent classes $U_i(\Lambda)$, $i\in\{1,\ldots,7\}$, which satisfy:
\begin{eqnarray*}
U_1(\Lambda)\cap X_{C_4}(\Lambda)\hspace{-0.2cm} &= \hspace{-0.2cm}&\{\pm{\bf v}_1=\pm(0,2,0,0),\ \pm{\bf v}_2=\pm(0,0,0,2)\},\\
U_2(\Lambda)\cap X_{C_4}(\Lambda)\hspace{-0.2cm} &= \hspace{-0.2cm}&\{\pm\mbox{$\frac{1}{2}$}({\bf v}_1+{\bf v}_2)=\pm(0,1,0,1),\ \pm{\bf v}_3=\pm(2,0,0,0)\},\\
U_3(\Lambda)\cap X_{C_4}(\Lambda)\hspace{-0.2cm} &= \hspace{-0.2cm}&\{\pm\mbox{$\frac{1}{2}$}({\bf v}_1-{\bf v}_2)=\pm(0,1,0,-1)\},\\
U_4(\Lambda)\cap X_{C_4}(\Lambda)\hspace{-0.2cm} &= \hspace{-0.2cm}&\{\pm\mbox{$\frac{1}{2}(\frac{1}{2}({\bf v}_1+{\bf v}_2)+{\bf v}_3)=\pm(1,\frac{1}{2},0,\frac{1}{2}$}),\ \pm{\bf v}_4=\pm(0,0,2,0)\},\\
U_5(\Lambda)\cap X_{C_4}(\Lambda)\hspace{-0.2cm} &= \hspace{-0.2cm}&\{\pm\mbox{$\frac{1}{2}(\frac{1}{2}({\bf v}_1+{\bf v}_2)-{\bf v}_3)=\pm(-1,\frac{1}{2},0,\frac{1}{2}$})\},\\
U_6(\Lambda)\cap X_{C_4}(\Lambda)\hspace{-0.2cm} &= \hspace{-0.2cm}&\{\pm\mbox{$\frac{1}{2}(\frac{1}{2}(\frac{1}{2}({\bf v}_1+{\bf v}_2)+{\bf v}_3)+{\bf v}_4)=\pm(\frac{1}{2},\frac{1}{4},1,\frac{1}{4}$})\},\\
U_7(\Lambda)\cap X_{C_4}(\Lambda)\hspace{-0.2cm} &= \hspace{-0.2cm}&\{\pm\mbox{$\frac{1}{2}(\frac{1}{2}(\frac{1}{2}({\bf v}_1+{\bf v}_2)+{\bf v}_3)-{\bf v}_4)=\pm(\frac{1}{2},\frac{1}{4},-1,\frac{1}{4}$})\},
\end{eqnarray*}
up to permutation of coordinates and signs of coordinates.

\medskip
\noindent
{\bf Case 3.1.1.} Scenario (C) happens.

It is easy to see that scenario (C) cannot happen twice time and for the remained eight modulo 2 equivalent classes of $\Lambda$, scenario (A) or (B) cannot happen. Therefore, we have
$${\rm card}\{X_{C_4}(\Lambda)\}\leq20+(4+2+2)\times2+2+2=40.$$

\medskip
\noindent
{\bf Case 3.1.2.} Scenario (B) happens and scenario (C) are not happen.

In this case, we use the notations of scenario (B) which discussed before. If scenario (B) happened one time and for the remained eight modulo 2 equivalent classes of $\Lambda$, scenario (A) are not happen, we have
$${\rm card}\{X_{C_4}(\Lambda)\}\leq20+(4+2+2)\times2+2+2=40.$$
If scenario (B) happened twice time, or scenario (B) and scenario (A) happened simultaneously, we have $\pm(2,0,0,0),\ \pm(0,0,2,0)\in\Lambda$. By
\begin{eqnarray*}
(2,0,0,0)\hspace{-0.2cm} &\pm \hspace{-0.2cm}&(v_{31},0,v_{33},0),\quad (2,0,0,0)\pm(v_{41},0,v_{43},0),\\
(0,0,2,0)\hspace{-0.2cm} &\pm \hspace{-0.2cm}&(v_{31},0,v_{33},0),\quad (0,0,2,0)\pm(v_{41},0,v_{43},0)\in\Lambda,
\end{eqnarray*}
we have $\{\pm{\bf v}_3, \pm{\bf v}_4\}=\{\pm(1,0,1,0), \pm(1,0,-1,0)\}$, which means the lattice $\Lambda$ is generated by
$$(2,0,0,0),\quad (1,1,0,0),\quad (1,0,1,0),\quad (\mbox{$\frac{1}{2}, \frac{1}{2}, \frac{1}{2}, \frac{1}{2}$})$$
and
$${\rm card}\{X_{C_4}(\Lambda)\}=40.$$

\medskip
\noindent
{\bf Case 3.1.3.} Scenario (A) happens and scenario (B) and (C) are not happen.

If scenario (A) happens only one time, we have
$${\rm card}\{X_{C_4}(\Lambda)\}\leq14+(4+2+2)\times3+2=40.$$

If scenario (A) happens twice time and scenario (${\rm A_1}$) happens, we use the notations in scenario (${\rm A_1}$).

Suppose for the ten remained modulo 2 equivalent classes of $\Lambda$, scenario (${\rm A_1}$) happens, then there exist modulo 2 equivalent classes $U_6(\Lambda)$ and $U_7(\Lambda)$, which satisfy:
\begin{eqnarray*}
U_6(\Lambda)\cap X_{C_4}(\Lambda)\hspace{-0.2cm} &= \hspace{-0.2cm}&\{\pm{\bf v}_1'=\pm(0,0,2,0),\ \pm{\bf v}_2'=\pm(v_{21}',v_{22}',0,0)\},\\
U_7(\Lambda)\cap X_{C_4}(\Lambda)\hspace{-0.2cm} &= \hspace{-0.2cm}&\{\pm\mbox{$\frac{1}{2}$}({\bf v}_1'+{\bf v}_2')=\pm(\mbox{$\frac{1}{2}$}v_{21}',\mbox{$\frac{1}{2}$}v_{22}',1,0),\ \pm{\bf v}_3'=\pm(0,0,0,2)\},
\end{eqnarray*}
where $|v_{21}'|+|v_{22}'|=2$. By ${\bf v}_1\pm{\bf v}_2'$, ${\bf v}_3\pm{\bf v}_2'\in\Lambda$, we have $|v_{21}'|=|v_{22}'|=1$. Therefore, we have
$$\min\big\{||{\bf v}_2'\pm\mbox{$\frac{1}{2}(\frac{1}{2}$}({\bf v}_1+{\bf v}_2)+{\bf v}_3)||_{C_4},\ ||{\bf v}_2'\pm\mbox{$\frac{1}{2}(\frac{1}{2}$}({\bf v}_1+{\bf v}_2)-{\bf v}_3)||_{C_4}\big\}=1,$$
which is a contradiction.

Suppose for the ten remained modulo 2 equivalent classes of $\Lambda$, scenario (${\rm A_2}$) happens. Then there exist modulo 2 equivalent classes $U_6(\Lambda)$ and $U_7(\Lambda)$, which satisfy:
\begin{eqnarray*}
U_6(\Lambda)\cap X_{C_4}(\Lambda)\hspace{-0.2cm} &= \hspace{-0.2cm}&\{\pm{\bf v}_1'=\pm(0,0,2,0),\ \pm{\bf v}_2'=\pm(0,0,0,-2)\},\\
U_7(\Lambda)\cap X_{C_4}(\Lambda)\hspace{-0.2cm} &= \hspace{-0.2cm}&\{\pm\mbox{$\frac{1}{2}$}({\bf v}_1'+{\bf v}_2')=\pm(0,0,1,-1),\ \pm{\bf v}_3'=\pm(v_{31}',v_{32}',0,0)\},
\end{eqnarray*}
where $|v_{31}'|+|v_{32}'|=2$. By ${\bf v}_1\pm{\bf v}_3'$, ${\bf v}_3\pm{\bf v}_3'\in\Lambda$, we have $|v_{31}'|=|v_{32}'|=1$, which is also a contradiction.

If scenario (${\rm A_2}$) happens twice time, one can easily verify that, the lattice $\Lambda$ is generated by
$$(2,0,0,0),\quad (1,1,0,0),\quad (1,0,1,0),\quad (\mbox{$\frac{1}{2}, \frac{1}{2}, \frac{1}{2}, \frac{1}{2}$})$$
and
$${\rm card}\{X_{C_4}(\Lambda)\}=40.$$

\medskip
\noindent
{\bf Case 3.2.} There exist modulo 2 equivalent class of $\Lambda$ can contain three pairs vectors $\pm{\bf v}_1', \pm{\bf v}_2', \pm{\bf v}_3'$ of $X_{C_4}(\Lambda)$.

By (3.1), we set ${\bf v}_1'=(2,0,0,0),\ {\bf v}_2'=(0,2,0,0),\ {\bf v}_3'=(0,0,v_{33}',v_{34}')$ with $v_{33}',\ v_{34}'\geq0$, $v_{33}'+v_{34}'=2$, without loss of generality. Then there exist modulo 2 equivalent classes $U_i(\Lambda)$, $i\in\{1,\ldots,7\}$, which satisfy:
\begin{align}
& U_1(\Lambda)\cap X_{C_4}(\Lambda)=\{\pm{\bf v}_1'=\pm(2,0,0,0),\ \pm{\bf v}_2'=\pm(0,2,0,0),\ \pm{\bf v}_3'=\pm(0,0,v_{33}',v_{34}')\},\notag\\
& U_2(\Lambda)\cap X_{C_4}(\Lambda)\supset\{\pm\mbox{$\frac{1}{2}$}({\bf v}_1'+{\bf v}_2')=\pm(1,1,0,0)\}, \tag{4.4}\\ 
& U_3(\Lambda)\cap X_{C_4}(\Lambda)\supset\{\pm\mbox{$\frac{1}{2}$}({\bf v}_1'-{\bf v}_2')=\pm(1,-1,0,0)\}, \tag{4.5}\\
& U_4(\Lambda)\cap X_{C_4}(\Lambda)\supset\{\pm\mbox{$\frac{1}{2}$}({\bf v}_1'+{\bf v}_3')=\pm(1,0,\mbox{$\frac{1}{2}$}v_{33}',\mbox{$\frac{1}{2}$}v_{34}')\}, \tag{4.6}\\
& U_5(\Lambda)\cap X_{C_4}(\Lambda)\supset\{\pm\mbox{$\frac{1}{2}$}({\bf v}_1'-{\bf v}_3')=\pm(1,0,-\mbox{$\frac{1}{2}$}v_{33}',-\mbox{$\frac{1}{2}$}v_{34}')\}, \tag{4.7}\\
& U_6(\Lambda)\cap X_{C_4}(\Lambda)\supset\{\pm\mbox{$\frac{1}{2}$}({\bf v}_2'+{\bf v}_3')=\pm(0,1,\mbox{$\frac{1}{2}$}v_{33}',\mbox{$\frac{1}{2}$}v_{34}')\}, \tag{4.8}\\
& U_7(\Lambda)\cap X_{C_4}(\Lambda)\supset\{\pm\mbox{$\frac{1}{2}$}({\bf v}_2'-{\bf v}_3')=\pm(0,1,-\mbox{$\frac{1}{2}$}v_{33}',-\mbox{$\frac{1}{2}$}v_{34}')\}. \tag{4.9}
\end{align}

If all the \lq\lq $\supset$" in (4.4)-(4.9) are \lq\lq =", then scenario (A) or (B) must happens in the remained eight modulo 2 equivalent classes of $\Lambda$. If scenario (A) happens, we have
$${\rm card}\{X_{C_4}(\Lambda)\}\leq18+14+4+2+2=40;$$
if scenario (B) happens, we have
$${\rm card}\{X_{C_4}(\Lambda)\}\leq18+20+2=40.$$

From now on, we assume one of \lq\lq $\supset$" in (4.4)-(4.9) cannot be \lq\lq =".

\medskip
\noindent
{\bf Case 3.2.1.} $v_{33}'\cdot v_{34}'\neq0$.

\medskip
By (3.1), both $\pm\frac{1}{2}({\bf v}_1'+{\bf v}_3')$, $\pm\frac{1}{2}({\bf v}_1'-{\bf v}_3')$, $\pm\frac{1}{2}({\bf v}_2'+{\bf v}_3')$ and $\pm\frac{1}{2}({\bf v}_2'-{\bf v}_3')$ cannot modulo 2 equivalent to other vectors (besides themselves, respectively) of $X_{C_4}(\Lambda)$. Therefore, we assume ${\bf v}_4'=(0,0,v_{43}',v_{44}')$ are modulo 2 equivalent to $\pm\mbox{$\frac{1}{2}$}({\bf v}_1'+{\bf v}_2')=\pm(1,1,0,0)$, by (3.1).

By the discussions of scenario (B), both the mid-point of $\mbox{$\frac{1}{2}$}({\bf v}_1'+{\bf v}_2')$ and ${\bf v}_4'$ and the mid-point of $(-\frac{1}{2}({\bf v}_1'+{\bf v}_2'))$ and ${\bf v}_4'$ cannot modulo 2 equivalent to other vectors of $X_{C_4}(\Lambda)$. In other words, there exist modulo 2 equivalent classes $U_8(\Lambda), U_9(\Lambda)$, which satisfy
\begin{eqnarray*}
U_8(\Lambda)\cap X_{C_4}(\Lambda)\hspace{-0.2cm} &= \hspace{-0.2cm}&\{\pm\mbox{$\frac{1}{2}(\frac{1}{2}({\bf v}_1'+{\bf v}_2')+{\bf v}_4')=\pm(\frac{1}{2},\frac{1}{2}$},\mbox{$\frac{1}{2}$}v_{43}', \mbox{$\frac{1}{2}$}v_{44}')\},\\
U_9(\Lambda)\cap X_{C_4}(\Lambda)\hspace{-0.2cm} &= \hspace{-0.2cm}&\{\pm\mbox{$\frac{1}{2}(\frac{1}{2}({\bf v}_1'+{\bf v}_2')-{\bf v}_4')=\pm(\frac{1}{2},\frac{1}{2}$},-\mbox{$\frac{1}{2}$}v_{43}', -\mbox{$\frac{1}{2}$}v_{44}')\}.
\end{eqnarray*}

\medskip
\noindent
{\bf Case 3.2.1.1.} The \lq\lq $\supset$" in (4.5) are \lq\lq =".

If scenario (A) happens in remained six modulo 2 equivalent classes of $\Lambda$, we have
$${\rm card}\{X_{C_4}(\Lambda)\}\leq24+14+2=40;$$
if scenario (A) are not happens in remained six modulo 2 equivalent classes of $\Lambda$, we have
$${\rm card}\{X_{C_4}(\Lambda)\}\leq24+(4+2+2)\times2=40.$$

\medskip
\noindent
{\bf Case 3.2.1.2.} There exist ${\bf v}_5'=(0,0,v_{53}',v_{54}')\in X_{C_4}(\Lambda)$ are modulo 2 equivalent to $\pm\frac{1}{2}({\bf v}_1'-{\bf v}_2')=\pm(1,-1,0,0)$.

Similar to the discussions before, there exist modulo 2 equivalent classes $U_{10}(\Lambda), U_{11}(\Lambda)$, which satisfy:
\begin{eqnarray*}
U_{10}(\Lambda)\cap X_{C_4}(\Lambda)\hspace{-0.2cm} &= \hspace{-0.2cm}&\{\pm\mbox{$\frac{1}{2}(\frac{1}{2}({\bf v}_1'-{\bf v}_2')+{\bf v}_5')=\pm(\frac{1}{2},-\frac{1}{2}$},\mbox{$\frac{1}{2}$}v_{53}',\mbox{$\frac{1}{2}$}v_{54}')\},\\
U_{11}(\Lambda)\cap X_{C_4}(\Lambda)\hspace{-0.2cm} &= \hspace{-0.2cm}&\{\pm\mbox{$\frac{1}{2}(\frac{1}{2}({\bf v}_1'-{\bf v}_2')-{\bf v}_5')=\pm(\frac{1}{2},-\frac{1}{2}$},-\mbox{$\frac{1}{2}$}v_{53}',-\mbox{$\frac{1}{2}$}v_{54}')\}
\end{eqnarray*}
and we have
$${\rm card}\{X_{C_4}(\Lambda)\}\leq22+8+(4+2+2)+2=40.$$

\medskip
\noindent
{\bf Case 3.2.2.} $v_{33}'=2,\ v_{34}'=0.$

It is easy to see that, for the remained eight modulo 2 equivalent classes of $\Lambda$, both scenario (A), (B) and (C) cannot happen.

If there are at most two of \lq\lq $\supset$" in (4.4)-(4.9) are not \lq\lq =", by the results of Case 3.2.1, we have
$${\rm card}\{X_{C_4}(\Lambda)\}\leq40.$$

If there are exactly three of \lq\lq $\supset$" in (4.4)-(4.9) are not \lq\lq =", we have
$${\rm card}\{X_{C_4}(\Lambda)\}\leq24+6\times2+2+2=40.$$

If there are at least four of \lq\lq $\supset$" in (4.4)-(4.9) are not \lq\lq =", by (3.1), without loss of generality, we may assume there exist ${\bf v}_4'=(0,0,v_{43}',v_{44}'),\ {\bf v}_5'=(0,v_{52}',0,v_{54}')\in X_{C_4}(\Lambda)$ with $|v_{43}'|+|v_{44}'|=|v_{52}'|+|v_{54}'|=2$, which are modulo 2 equivalent to $\pm(1,1,0,0),\ \pm(1,0,1,0)$, respectively. Therefore, we have
$$\pm\mbox{$\frac{1}{2}(\frac{1}{2}$}({\bf v}_1'+{\bf v}_2')+{\bf v}_4')=\pm(\mbox{$\frac{1}{2},\frac{1}{2}$},\mbox{$\frac{1}{2}$}v_{43}',\mbox{$\frac{1}{2}$}v_{44}'),\quad \pm\mbox{$\frac{1}{2}(\frac{1}{2}$}({\bf v}_1'+{\bf v}_2')-{\bf v}_4')=\pm(\mbox{$\frac{1}{2},\frac{1}{2}$},-\mbox{$\frac{1}{2}$}v_{43}',-\mbox{$\frac{1}{2}$}v_{44}'),$$
$$\pm\mbox{$\frac{1}{2}(\frac{1}{2}$}({\bf v}_1'+{\bf v}_3')+{\bf v}_5')=\pm(\mbox{$\frac{1}{2},\mbox{$\frac{1}{2}$}v_{52}',\frac{1}{2}$},\mbox{$\frac{1}{2}$}v_{54}'),\quad \pm\mbox{$\frac{1}{2}(\frac{1}{2}$}({\bf v}_1'+{\bf v}_3')-{\bf v}_5')=\pm(\mbox{$\frac{1}{2},-\mbox{$\frac{1}{2}$}v_{52}',\frac{1}{2}$},-\mbox{$\frac{1}{2}$}v_{54}')$$
are belongs to $\Lambda$. By ${\bf v}_2'\pm{\bf v}_5'$, ${\bf v}_3'\pm{\bf v}_4'\in\Lambda$, we have $|v_{52}'|,\ |v_{43}'|\leq1$.

If $v_{52}'\geq0$, combining with
$$\mbox{$\frac{1}{2}(\frac{1}{2}$}({\bf v}_1'+{\bf v}_2')+{\bf v}_4')-\mbox{$\frac{1}{2}(\frac{1}{2}$}({\bf v}_1'+{\bf v}_3')+{\bf v}_5'),\ \mbox{$\frac{1}{2}(\frac{1}{2}$}({\bf v}_1'+{\bf v}_2')-{\bf v}_4')-\mbox{$\frac{1}{2}(\frac{1}{2}$}({\bf v}_1'+{\bf v}_3')+{\bf v}_5')\in\Lambda,$$
we have
$$\mbox{$\frac{1}{2}-\mbox{$\frac{1}{2}$}v_{52}'+\frac{1}{2}$}-\mbox{$\frac{1}{2}$}v_{43}'+|\mbox{$\frac{1}{2}$}v_{44}'-\mbox{$\frac{1}{2}$}v_{54}'|\geq2,\quad \mbox{$\frac{1}{2}-\mbox{$\frac{1}{2}$}v_{52}'+\frac{1}{2}$}+\mbox{$\frac{1}{2}$}v_{43}'+|\mbox{$\frac{1}{2}$}v_{44}'+\mbox{$\frac{1}{2}$}v_{54}'|\geq2,$$
which means
$$-v_{52}'+|v_{44}'|\geq2\quad {\rm or}\quad -v_{52}'+|v_{54}'|\geq2.$$
Therefore, we have $v_{52}'=0,\ \pm{\bf v}_5'=\pm(0,0,0,2)$.

By ${\bf v}_4'\pm{\bf v}_5'\in\Lambda$, we have $\pm{\bf v}_4'=\pm(0,0,1,1)$ or $\pm(0,0,1,-1)$, the mid-point of ${\bf v}_3'$ and ${\bf v}_5'$ or the mid-point of ${\bf v}_3'$ and $(-{\bf v}_5')$. In other words, there exist a modulo 2 equivalent class of $\Lambda$ which can contain four pairs vectors of $X_{C_4}(\Lambda)$, which is a contradiction.

By consider $\mbox{$\frac{1}{2}(\frac{1}{2}$}({\bf v}_1'+{\bf v}_3')-{\bf v}_5')$ instead of $\mbox{$\frac{1}{2}(\frac{1}{2}$}({\bf v}_1'+{\bf v}_3')+{\bf v}_5')$, the $v_{52}'<0$ case also leads to a contradiction.

\medskip
As a conclusion of all the cases, we have
$$\kappa^*(C_4)=40,$$
Theorem 3 is proved.

\end{proof}

\medskip
Denote $\Lambda$ be the lattice which generated by $\{(-\frac{3}{4},\frac{3}{4},\frac{1}{4},\frac{1}{4}),\ (\frac{1}{4},-\frac{3}{4},\frac{3}{4},\frac{1}{4}),\ (\frac{1}{4},\frac{1}{4},-\frac{3}{4},\frac{3}{4}),$\\ $(\frac{3}{4},\frac{1}{4},\frac{1}{4},-\frac{3}{4})\}.$ One can easily verify that $C_4+\Lambda$ is a packing, ${\rm card}\{X_{C_4}(\Lambda)\}=26$ and the packing density of $C_4+\Lambda$ is $32/45$. Therefore, we have the following Remark.

\medskip
\noindent
{\bf Remark 2.} Similar to the tetrahedral case (see \cite{Zo02, Zo03}),  the density of the lattice packing $C_4+\Lambda$ attaining $\kappa^*(C_4)=40$ is only $2/3$ while the density of some lattice packing $C_4+\Lambda$ in which every one touches 26 others is $32/45$.

\medskip
We end this section by two conjectures as following.

\medskip
\noindent
{\bf Conjecture 4.1.} {\it In $\mathbb{E}^4$, we have
$$\delta^*(C_4)=32/45,$$
and it can be attained if and only if the corresponding lattice is generated by $\{(-\frac{3}{4},\frac{3}{4},\frac{1}{4},\frac{1}{4}),$ $(\frac{1}{4},-\frac{3}{4},\frac{3}{4},\frac{1}{4}), (\frac{1}{4},\frac{1}{4},-\frac{3}{4},\frac{3}{4}), (\frac{3}{4},\frac{1}{4},\frac{1}{4},-\frac{3}{4})\},$ up to permutation of coordinates and signs of coordinates.}

\medskip
\noindent
{\bf Conjecture 4.2.} {\it In $\mathbb{E}^4$, we have
$$\kappa(C_4)=40,$$
and it can be attained if and only if by the lattice $\Lambda$ which generated by $\{(2,0,0,0),$ $(1,1,0,0), (1,0,1,0), (\frac{1}{2},\frac{1}{2},\frac{1}{2},\frac{1}{2})\}$, up to permutation of coordinates and signs of coordinates.}

\vspace{1cm}
\noindent
{\LARGE\bf 5. The lattice kissing numbers of tetrahedra}

\bigskip
\noindent
Let $T$ denote the regular tetrahedron with vertices $(-\frac{1}{2}, \frac{1}{2}, \frac{1}{2}),\ (\frac{1}{2}, -\frac{1}{2}, \frac{1}{2}),\ (\frac{1}{2}, \frac{1}{2}, -\frac{1}{2})$ and $(-\frac{1}{2}, -\frac{1}{2}, -\frac{1}{2})$. In 1996, Zong \cite{Zo02} proved that
$$\kappa^*(T)=18\eqno(5.1)$$
and
$$\kappa(T)=18\quad {\rm or}\quad 19.$$
Later, in 1999, by modifying Zong's method, Talata \cite{Ta01} proved
$$\kappa(T)=18.\eqno(5.2)$$

In 1904, Minkowski \cite{Mi01} defined
$$D(K)=\{{\bf x}-{\bf y}:\ {\bf x}-{\bf y}\in K\}$$
for a convex body $K$. Usually, we call $D(K)$ the different set of $K$.

Clearly, $D(K)$ is a centrally symmetric convex body centered at ${\bf o}$. Furthermore, Minkowski showed that
$$(K+{\bf x})\cap(K+{\bf y})\neq\emptyset$$
if and only if
$$(\mbox{$\frac{1}{2}$}D(K)+{\bf x})\cap(\mbox{$\frac{1}{2}$}D(K)+{\bf y})\neq\emptyset.$$
Therefore, one can easily deduce that
$$\kappa^*(K)=\kappa^*(D(K))$$
and
$$\kappa(K)=\kappa(D(K)).$$

It is well known that
$$D(T)=\{(v_1, v_2, v_3):\ |v_1|+|v_2|+|v_3|\leq2,\ |v_1|, |v_2|, |v_3|\leq1\}.$$
Usually, we call $D(T)$ a cubeoctahedron.

In this section, we give a rather simpler proof of (5.1).

\medskip
\noindent
{\bf Theorem 5.1.} {\it For regular tetrahedron $T$, we have
$$\kappa^*(T)=18.$$}
\begin{proof}
It is easy to see that
$$B^3\subset D(T)\subset\sqrt2B^3.$$
Let $D(T)+\Lambda$ be a lattice packing which attaining $\kappa^*(D(T))$ and
$$X=\{{\bf v}_1,...,{\bf v}_{\kappa^*(D(T))}\}=\partial(2D(T))\cap\Lambda,$$
where $\partial(2D(T))$ denotes the boundary of $2D(T)$. Since
$$2B^3\subset 2D(T)\subset 2\sqrt2B^3,$$
we have
$$2\leq||{\bf v}_i||\leq2\sqrt2$$
holds for $i=1,2,...,\kappa^*(D(T))$. Since $\Lambda$ is also a packing lattice of $B^3$, by \cite{Li01} we get
$$\kappa^*(D(T))={\rm card}X\leq20,$$
and ${\rm card}X=20$ if and only if $\Lambda$ is generated by
$${\bf a}_1=(2,0,0),\quad {\bf a}_2=(0,2,0)\quad {\rm and}\quad {\bf a}_3=(1,0,\sqrt3),$$
up to rotation and reflection; ${\rm card}X=18$ if and only if $\Lambda$ is generated by
$${\bf a}_1=(2, 0, 0),\quad {\bf a}_2=(0, 2, 0)\quad {\rm and}\quad {\bf a}_3=(1, 1, \sqrt2)$$
or
$${\bf a}_1=(2, 0, 0),\quad {\bf a}_2=(0, 2, 0)\quad {\rm and}\quad {\bf a}_3=(0, 0, 2),$$
up to rotation and reflection.
By verify three lattices above, one can easily check that
$$\kappa^*(D(T))=18,$$
and it can be attained if and only if the corresponding lattice is generated by ${\bf a}_1=(2, 0, 0),\ {\bf a}_2=(0, 2, 0)\ {\rm and}\ {\bf a}_3=(0, 0, 2).$ Therefore, we have
$$\kappa^*(T)=\kappa^*(D(T))=18,$$
Theorem 5.1 is proved.
\end{proof}

\vspace{1cm}
\noindent
{\LARGE\bf 6. A Link Between Kissing Numbers and Packing Densities}

\bigskip
\noindent
As we mentioned in Remark 2, for a centrally symmetric convex body $C$, the lattices which attaining $\kappa^*(C)$ and $\delta^*(C)$ sometimes are not identity. Thus, the kissing number and the packing density of $C$ seems cannot determined simultaneously. However, one may imagine that, there still exist some relations between them.

By generalizing Blichfeldt's method \cite{Bl01}, Rankin \cite{Ra01} proved the following results:
$$\kappa(B^n)\ll\frac{\sqrt{\pi n^3}}{2}2^{0.5n} \eqno(6.1)$$
and
$$\delta(B^n)\ll\frac{n+2}{2}2^{-0.5n}, \eqno(6.2)$$
where $\delta(B^n)$ is the translative packing density of $B^n$ and $f(x)\ll g(x)$ means \\$\limsup\limits_{x\to\infty}\frac{f(x)}{g(x)}\leq 1.$
In 1978, Kabatjanski and Levenshtein \cite{Ka01} studied the Jacobi polynomials and spherical polynomials deeply, combining with Delsarte's Lemma (see \cite{De01}) and lots of real analysis methods, they proved that
$$\kappa(B^n)\leq2^{0.401n(1+o(1))} \eqno(6.3)$$
and
$$\delta(B^n)\leq2^{-0.599n(1+o(1))}. \eqno(6.4)$$

By (6.1)-(6.4), one may conjecture that $\delta(B^n)\ll\kappa(B^n)/2^n.$ However, until now, this relationship is unproven. For more details and references concerning the upper bounds for $\delta(B^n)$ and $\kappa(B^n)$, we refer to \cite{Bo01, Co01, Co02, Fe01, Ro01, Ro02, Sc03, Zo01}.

For an $n$-dimensional centrally symmetric convex body $C$ and a positive number $\rho$, define
$$\tau_\rho(C,X)=\max_{{\bf x}\in\mathbb{E}^n}{\rm card}\{{\bf v}:\ ||{\bf v}, {\bf x}||_C\leq\rho,\ {\bf v}\in X\}$$
and
$$\tau_\rho(C)=\max_{X}\tau_\rho(C,X),$$
where $X$ is a discrete set such that $C+X$ is a packing in $\mathbb{E}^n$. For $\rho\geq2$, it is easy to see that
$$\tau_\rho(C)\geq\kappa_\alpha(C)+1,\quad \alpha=\rho-2,$$
where $\kappa_\alpha(C)$ is the generalized kissing numbers of $C$, namely: the maximum number of translative sets of $C$ that can be packed in $(3+\alpha)C\setminus {\rm int}(C)$, see \cite{Li01}.

\medskip
\noindent
{\bf Lemma 6.1.} {\it For all centrally symmetric convex body $C$ centered at {\bf o} and ${\bf v}\in\partial C$, we have
$$C\setminus({\rm int}C+\mu{\bf v})\subset C\setminus({\rm int}C+{\bf v})$$
holds for $0<\mu\leq1$.}
\begin{proof}
For a point ${\bf x}\in C\cap({\rm int}C+{\bf v})$, we have
\begin{eqnarray*}
||{\bf x}-\mu{\bf v}||_C\hspace{-0.2cm}&=\hspace{-0.2cm}&||(1-\mu){\bf x}+\mu({\bf x}-{\bf v})||_C\\
\hspace{-0.2cm}&\leq\hspace{-0.2cm}& (1-\mu)||{\bf x}||_C+\mu||{\bf x}-{\bf v}||_C\\
\hspace{-0.2cm}&<\hspace{-0.2cm}& 1.
\end{eqnarray*}
Therefore, we have
$${\bf x}\in C\cap({\rm int}C+\mu{\bf v}),$$
which means that
$$C\cap({\rm int}C+\mu{\bf v})\supset C\cap({\rm int}C+{\bf v}).$$

In other words,
$$C\setminus({\rm int}C+\mu{\bf v})\subset C\setminus({\rm int}C+{\bf v}),$$
Lemma 6.1 is proved.
\end{proof}

\medskip
\noindent
{\bf Lemma 6.2.} {\it For all centrally symmetric convex body $C$, we have
$$\tau_2(C)=\kappa(C)+1.$$}
\begin{proof}
By the definition of $\kappa(C)$, we have
$$\kappa(C)=\max_{X}{\rm card}\{{\bf v}:\ ||{\bf v}||_C=2,\ {\bf v}\in X\}, \eqno(6.5)$$
where $X$ is a discrete set such that $C+X$ is a packing in $\mathbb{E}^n$.

For a packing set $X_0$ of $C$, we suppose
$$\tau_2(C, X_0)={\rm card}\{{\bf v}:\ ||{\bf v}||_C\leq2,\ {\bf v}\in X_0\},$$
without loss of generality.

If ${\bf o}\in X_0$, one can easily obtain
$$\tau_2(C, X_0)\leq\kappa(C)+1,$$
and the equality holds if and only if there exist a finite subset $X_0'\subset X_0$ such that $C+X_0'$ attaining the translative kissing number $\kappa(C)$ of $C$.

On the contrary, if ${\bf o}\notin X_0$, we denote
$$\{{\bf v}:\ ||{\bf v}||_C\leq2,\ {\bf v}\in X_0\}=\{{\bf v}_1, {\bf v}_2,\ldots, {\bf v}_{\tau_2(C, X_0)}\},$$
namely
$$0<||{\bf v}_1||_C, ||{\bf v}_2||_C, \ldots,||{\bf v}_{\tau_2(C, X_0)}||_C\leq2 \eqno(6.6)$$
and
$$||{\bf v}_i, {\bf v}_j||_C\geq2 \eqno(6.7)$$
holds for $i\neq j$. Denote ${\bf v}_i'=\frac{2}{||{\bf v}_i||_C}{\bf v}_i$ for all $i$, we have $||{\bf v}_i'||_C=2.$ Combining with (6.6), (6.7) and Lemma 6.1, we have
$${\bf v}_j\in 2C\setminus\big(\{{\bf o}\}\cup({\rm int}2C+{\bf v}_1)\big)\subset2C\setminus\big(\{{\bf o}\}\cup({\rm int}2C+{\bf v}_1')\big)$$
holds for all $j\neq1$. Therefore, let ${\bf v}_1'$ substitute ${\bf v}_1$, (6.6) and (6.7) still holds. By inductive, we have
$$||{\bf v}_1'||_C=||{\bf v}_2'||_C=\ldots=||{\bf v}_{\tau_2(C, X_0)}'||_C=2 \eqno(6.6')$$
and
$$||{\bf v}_i', {\bf v}_j'||_C\geq2 \eqno(6.7')$$
holds for $i\neq j$. By (6.5), we have
$$\tau_2(C, X_0)\leq\kappa(C).$$

As a conclusion of the two cases, we have
$$\tau_2(C)=\kappa(C)+1,$$
Lemma 6.2 is proved.
\end{proof}

\medskip
\noindent
{\bf Lemma 6.3. (Blichfeldt \cite{Bl01})} {\it Let $\varphi\left({\bf v}\right)$ be a bounded, nonnegative Riemann integrable function and vanishing outside a bounded region. Suppose that for each translative packing $C+X$,
$$\sum_{{\bf v}_i\in X}\varphi\left({\bf v}-{\bf v}_i\right)\leq 1 \eqno(6.8)$$
holds for every point ${\bf v}\in \mathbb{E}^{n}$. Then
$$\delta(C)\leq {\rm vol}\left(C\right)\times {1\over \int_{\mathbb{R}^n}\varphi\left({\bf v}\right)d{\bf v}}.$$}

\medskip
For an $n$-dimensional centrally symmetric convex body $C$ and a positive number $\rho$, define
$$\varphi\left({\bf v}\right)=\left\{
\begin{aligned}
\frac{1}{\tau_\rho(C)}, \ \ {\bf v}\in\rho C,\\
0,\ \ {\bf v}\notin\rho C.\\
\end{aligned}
\right.
$$
It is easy to see that, this function $\varphi\left({\bf v}\right)$ satisfying (6.8). By Lemma 6.3, we have
$$\delta(C)\leq\frac{\tau_\rho(C)}{{\rho}^n}.$$
Combining with Lemma 6.2, we have the following theorem.

\medskip
\noindent
{\bf Theorem 4.} {\it For all $n$-dimensional centrally symmetric convex body $C$, we have
$$\delta(C)\leq\frac{\kappa(C)+1}{2^n}.$$}

\medskip
For an $n$-dimensional centrally symmetric convex body $C$ and a positive number $\rho$, define
$$\tau^*_\rho(C,\Lambda)=\max_{{\bf x}\in\mathbb{E}^n}{\rm card}\{{\bf v}:\ ||{\bf v}, {\bf x}||_C\leq\rho,\ {\bf v}\in \Lambda\}$$
and
$$\tau^*_\rho(C)=\max_{\Lambda}\tau^*_\rho(C,\Lambda),$$
where $\Lambda$ is a packing lattice of $C$. It is easy to see that 
$$\tau^*_\rho(C)\leq\tau_\rho(C).$$

 Similar to the translative case, for $\rho\geq2$, we have
$$\tau^*_\rho(C)\geq\kappa^*_\alpha(C)+1,\quad \alpha=\rho-2$$
and
$$\delta^*(C)\leq\frac{\tau^*_\rho(C)}{\rho^n}.$$

\medskip
In \cite{Micc02} (Chapter 5), Micciancio proved the following results: {\it If $\rho<\sqrt2$, then $\tau^*_\rho(B^n)$ is a constant independent of $n$; if $\rho=\sqrt2$, then $\tau^*_\rho(B^n)=2n$; if $\rho>\sqrt2$, then $\tau^*_\rho(B^n)>2^{n^c}$ for some fixed constant $c$ independent of $n$.}

\medskip
We end this section by a problem and a remark as following.

\medskip
\noindent
{\bf Problem 6.1.} {\it Is it true that
$$\tau^*_2(B^n)=\kappa^*(B^n)+1$$
holds for $n$-dimensional unit ball $B^n$?}

\medskip
\noindent
{\bf Remark 6.1.} If the answer of Problem 6.1 is positive, then we have
$$\delta^*(B^n)\leq\frac{\kappa^*(B^n)+1}{2^n}$$
and
$$\kappa^*(B^n)>2^{n^c}$$
holds for $n$-dimensional unit ball $B^n$, where $c$ is a constant independent of $n$. 

On the other hand, if the answer of Problem 6.1 is negative, the counter-example itself will be interesting.

\vspace{0.8cm}\noindent
{\bf Acknowledgement.} This work is supported by the National Natural Science Foundation of China (NSFC12226006, NSFC11921001) and the Natural Key Research and Development Program of China (2018YFA0704701).

\vspace{0.6cm}
\noindent
Yiming Li, Center for Applied Mathematics, Tianjin University, Tianjin 300072, China
\noindent
Email: xiaozhuang@tju.edu.cn

\medskip\noindent
Chuanming Zong, Center for Applied Mathematics, Tianjin University, Tianjin 300072, China

\noindent
Email: cmzong@tju.edu.cn

\end{document}